\documentclass[12pt]{article}
\usepackage{amsthm, amsmath, amssymb, enumerate}
\usepackage{xcolor}
\usepackage{fullpage}
\usepackage{stfloats}
\usepackage[colorlinks=true]{hyperref}
\usepackage{tikz}
\usetikzlibrary{positioning}
\usepackage[all]{xy}
\usepackage{esint}

\begin{document}
\title{Modular Heights of Unitary Shimura Varieties I: Derivatives of Eisenstein Series}
\author{Ziqi Guo}
\maketitle

\theoremstyle{plain}
\newtheorem{thm}{Theorem}[section]
\newtheorem{theorem}[thm]{Theorem}
\newtheorem{cor}[thm]{Corollary}
\newtheorem{corollary}[thm]{Corollary}
\newtheorem{lem}[thm]{Lemma}
\newtheorem{lemma}[thm]{Lemma}
\newtheorem{pro}[thm]{Proposition}
\newtheorem{proposition}[thm]{Proposition}
\newtheorem{prop}[thm]{Proposition}
\newtheorem{definition}[thm]{Definition}
\newtheorem{assumption}[thm]{Assumption}
\def\avint{\mathop{\,\rlap{-}\!\!\int}\nolimits}

\theoremstyle{remark}
\newtheorem{remark}[thm]{Remark}
\newtheorem{example}[thm]{Example}
\newtheorem{remarks}[thm]{Remarks}
\newtheorem{problem}[thm]{Problem}
\newtheorem{exercise}[thm]{Exercise}
\newtheorem{situation}[thm]{Situation}
\newtheorem{acknowledgment}[thm]{Acknowledgment}
\newtheorem*{remark*}{Remark}

\numberwithin{equation}{subsection}

\newcommand{\ZZ}{\mathbb{Z}}
\newcommand{\CC}{\mathbb{C}}
\newcommand{\QQ}{\mathbb{Q}}
\newcommand{\RR}{\mathbb{R}}
\newcommand{\HH}{\mathcal{H}}     

\newcommand{\ad}{\mathrm{ad}}            
\newcommand{\NT}{\mathrm{NT}}
\newcommand{\nonsplit}{\mathrm{nonsplit}}
\newcommand{\Pet}{\mathrm{Pet}}
\newcommand{\Fal}{\mathrm{Fal}}
\newcommand{\Af}{\mathbb{A}_f}

\newcommand{\cs}{{\mathrm{cs}}}

\newcommand{\XU}{X_U}
\newcommand{\Fn}{F_v}
\newcommand{\LU}{L_U}
\newcommand{\LL}{\overline{\mathcal{L}}}
\newcommand{\OF}{\mathcal{O}_F}
\renewcommand{\OE}{\mathcal{O}_E}
\newcommand{\XXU}{\mathcal{X}_U}
\newcommand{\OA}{\underline{\Omega}_\mathcal{A}}
\newcommand{\OU}{\Omega_{\mathcal{X}_U/\mathbb{Z}[\frac{1}{n}]}}
\newcommand{\WA}{\underline{\omega}_\mathcal{A}}
\newcommand{\WU}{\omega_{\mathcal{X}_U/\mathbb{Z}[\frac{1}{n}]}}
\newcommand{\HHom}{\mathcal{H}\mathrm{om}}

\newcommand{\pair}[1]{\langle {#1} \rangle}
\newcommand{\wpair}[1]{\left\{{#1}\right\}}
\newcommand{\wh}{\widehat}
\newcommand{\wt}{\widetilde}

\newcommand\Spf{\mathrm{Spf}}

\newcommand{\lra}{{\longrightarrow}}

\newcommand{\matrixx}[4]
{\left( \begin{array}{cc}
  #1 &  #2  \\
  #3 &  #4  \\
 \end{array}\right)}        


\newcommand{\BA}{{\mathbb {A}}}
\newcommand{\BB}{{\mathbb {B}}}
\newcommand{\BC}{{\mathbb {C}}}
\newcommand{\BD}{{\mathbb {D}}}
\newcommand{\BE}{{\mathbb {E}}}
\newcommand{\BF}{{\mathbb {F}}}
\newcommand{\BG}{{\mathbb {G}}}
\newcommand{\BH}{{\mathbb {H}}}
\newcommand{\BI}{{\mathbb {I}}}
\newcommand{\BJ}{{\mathbb {J}}}
\newcommand{\BK}{{\mathbb {K}}}
\newcommand{\BL}{{\mathbb {L}}}
\newcommand{\BM}{{\mathbb {M}}}
\newcommand{\BN}{{\mathbb {N}}}
\newcommand{\BO}{{\mathbb {O}}}
\newcommand{\BP}{{\mathbb {P}}}
\newcommand{\BQ}{{\mathbb {Q}}}
\newcommand{\BR}{{\mathbb {R}}}
\newcommand{\BS}{{\mathbb {S}}}
\newcommand{\BT}{{\mathbb {T}}}
\newcommand{\BU}{{\mathbb {U}}}
\newcommand{\BV}{{\mathbb {V}}}
\newcommand{\BW}{{\mathbb {W}}}
\newcommand{\BX}{{\mathbb {X}}}
\newcommand{\BY}{{\mathbb {Y}}}
\newcommand{\BZ}{{\mathbb {Z}}}

\newcommand{\CA}{{\mathcal {A}}}
\newcommand{\CB}{{\mathcal {B}}}
\newcommand{\CD}{{\mathcal{D}}}
\newcommand{\CE}{{\mathcal {E}}}
\newcommand{\CF}{{\mathcal {F}}}
\newcommand{\CG}{{\mathcal {G}}}
\newcommand{\CH}{{\mathcal {H}}}
\newcommand{\CI}{{\mathcal {I}}}
\newcommand{\CJ}{{\mathcal {J}}}
\newcommand{\CK}{{\mathcal {K}}}
\newcommand{\CL}{{\mathcal {L}}}
\newcommand{\CM}{{\mathcal {M}}}
\newcommand{\CN}{{\mathcal {N}}}
\newcommand{\CO}{{\mathcal {O}}}
\newcommand{\CP}{{\mathcal {P}}}
\newcommand{\CQ}{{\mathcal {Q}}}
\newcommand{\CR }{{\mathcal {R}}}
\newcommand{\CS}{{\mathcal {S}}}
\newcommand{\CT}{{\mathcal {T}}}
\newcommand{\CU}{{\mathcal {U}}}
\newcommand{\CV}{{\mathcal {V}}}
\newcommand{\CW}{{\mathcal {W}}}
\newcommand{\CX}{{\mathcal {X}}}
\newcommand{\CY}{{\mathcal {Y}}}
\newcommand{\CZ}{{\mathcal {Z}}}

\newcommand{\ab}{{\mathrm{ab}}}
\newcommand{\Ad}{{\mathrm{Ad}}}
\newcommand{\an}{{\mathrm{an}}}
\newcommand{\Aut}{{\mathrm{Aut}}}

\newcommand{\Br}{{\mathrm{Br}}}
\newcommand{\bs}{\backslash}
\newcommand{\bbs}{\|\cdot\|}

\newcommand{\Ch}{{\mathrm{Ch}}}
\newcommand{\cod}{{\mathrm{cod}}}
\newcommand{\cont}{{\mathrm{cont}}}
\newcommand{\cl}{{\mathrm{cl}}}
\newcommand{\criso}{{\mathrm{criso}}}
\newcommand{\de}{{\mathrm{d}}}
\newcommand{\dR}{{\mathrm{dR}}}
\newcommand{\df}{\mathrm{det}^*}
\newcommand{\disc}{{\mathrm{disc}}}
\newcommand{\Div}{{\mathrm{Div}}}
\renewcommand{\div}{{\mathrm{div}}}
\newcommand{\Dh}{\widehat{\mathcal{D}}}
\newcommand{\Ei}{\mathrm{Ei}}
\newcommand{\Eis}{{\mathrm{Eis}}}
\newcommand{\End}{{\mathrm{End}}}

\newcommand{\Frob}{{\mathrm{Frob}}}

\newcommand{\Gal}{{\mathrm{Gal}}}
\newcommand{\GL}{{\mathrm{GL}}}
\newcommand{\GO}{{\mathrm{GO}}}
\newcommand{\GSO}{{\mathrm{GSO}}}
\newcommand{\GSp}{{\mathrm{GSp}}}
\newcommand{\GSpin}{{\mathrm{GSpin}}}
\newcommand{\GU}{{\mathrm{GU}}}
\newcommand{\BGU}{{\mathbb{GU}}}

\newcommand{\Has}{\mathrm{hasse}}
\newcommand{\Hom}{{\mathrm{Hom}}}
\newcommand{\Hol}{{\mathrm{Hol}}}
\newcommand{\HC}{{\mathrm{HC}}}
\newcommand{\id}{\mathrm{id}}
\newcommand{\Img}{{\mathrm{Im}}}
\newcommand{\Ind}{{\mathrm{Ind}}}
\newcommand{\ine}{\mathrm{ine}}
\newcommand{\inv}{{\mathrm{inv}}}
\newcommand{\Isom}{{\mathrm{Isom}}}

\newcommand{\Jac}{{\mathrm{Jac}}}
\newcommand{\JL}{{\mathrm{JL}}}

\newcommand{\Ker}{{\mathrm{Ker}}}
\newcommand{\KS}{{\mathrm{KS}}}

\newcommand{\Lie}{{\mathrm{Lie}}}
\renewcommand{\mod}{\mathrm{mod}}
\newcommand{\mm}{\mathfrak{m}}
\newcommand{\new}{{\mathrm{new}}}
\newcommand{\Nm}{\mathrm{Nm}}
\newcommand{\NS}{{\mathrm{NS}}}

\newcommand{\ord}{{\mathrm{ord}}}
\newcommand{\ol}{\overline}
\newcommand{\otf}{\otimes^*}
\newcommand{\rank}{{\mathrm{rank}}}

\newcommand{\PGL}{{\mathrm{PGL}}}
\newcommand{\PSL}{{\mathrm{PSL}}}
\newcommand{\Pic}{\mathrm{Pic}}
\newcommand{\Prep}{\mathrm{Prep}}
\newcommand{\Proj}{\mathrm{Proj}}
\renewcommand{\Pr}{\mathcal{P}r}
\newcommand{\Picc}{\mathcal{P}ic}

\newcommand{\ram}{\mathrm{ram}}
\renewcommand{\Re}{{\mathrm{Re}}}
\newcommand{\Res}{{\mathrm{Res}}}
\newcommand{\red}{{\mathrm{red}}}
\newcommand{\reg}{{\mathrm{reg}}}
\newcommand{\sm}{{\mathrm{sm}}}
\newcommand{\sing}{{\mathrm{sing}}}
\newcommand{\SL}{\mathrm{SL}}
\newcommand{\SLL}{\widetilde{\mathrm{SL}}}
\newcommand{\SO}{\mathrm{SO}}
\newcommand{\Sp}{\mathrm{Sp}}
\newcommand{\spl}{\mathrm{spl}}
\newcommand{\Sym}{{\mathrm{Sym}}}
\newcommand{\Spec}{\mathrm{Spec}}
\renewcommand{\ss}{\mathrm{ss}}
\newcommand{\tor}{{\mathrm{tor}}}
\newcommand{\tr}{{\mathrm{tr}}}

\newcommand{\ur}{{\mathrm{ur}}}
\newcommand{\U}{\mathrm{U}}
\newcommand{\UU}{\mathrm{U}(1,1)}
\newcommand{\vol}{{\mathrm{vol}}}

\newcommand{\ds}{\displaystyle}

\begin{abstract}
   This is the first of a series of three papers, in which we prove a formula expressing the modular height of a unitary Shimura variety over a CM number field in terms of the logarithmic derivative of the Hecke L-function associated with the CM extension. The main idea of our proof is to compare the holomorphic projection of the derivative of a certain mixed Eisenstein-theta series and the arithmetic degree of a generating series of divisors on unitary Shimura varieties.

   In this paper, we compute an explicit expression of the holomorphic projection of the derivative of a certain mixed Eisenstein-theta series.
\end{abstract}

\tableofcontents

\section{Introduction}\label{introduction}
The goal of this series of three papers (\cite{Guo1}, \cite{Guo2} and the current one) is to prove a formula expressing the modular height of a unitary Shimura variety over a CM field in terms of the logarithmic derivative of the Hecke L-function associated with the CM extension. Our work can be viewed as an extension of X. Yuan's work \cite{Yuan1}, which is based on the work Yuan–Zhang–Zhang \cite{YZZ2} on the Gross–Zagier formula, and the work Yuan–Zhang \cite{YZ1} on the averaged Colmez conjecture. All these works are in turn inspired by the pioneering work Gross–Zagier \cite{GZ} and some philosophies of Kudla’s program \cite{Kud1,Kud2,Kud3,Kud4}. This series of works all aim to calculate the arithmetic invariants of Shimura varieties using special values of L-functions.

In our work, we will study the generating series of divisors on unitary Shimura varieties and their arithmetic versions, comparing them with the derivative of mixed theta-Eisenstein series. Through a series of specific and intricate calculations, we will provide a precise formula for the modular height. For a complete introduction to our work, we refer to the introduction to the third paper \cite{Guo2} in this series.

This paper is the first one in this series. Its goal is to provide ingredients on the ``analytic side" needed for our proof. More precisely, we give the explicit definition of the so-called mixed theta–Eisenstein series, and, after fixing a suitable Schwartz function, derive an explicit formula for the holomorphic projection of its derivative, which we call the derivative series. At the same time, we will also introduce some technical notation and a key lemma which are very useful in the treatment of the derivatives of theta series and Eisenstein series, and which will play an important role in the final proof of the formula.

\subsection{Pseudo-theta series and Pseudo-Eisenstein series}
The notion of pseudo-theta series and pseudo-Eisenstein series is learned from \cite{Yuan1}, and can be traced back to \cite{YZZ2} and \cite{YZ1}. Roughly speaking, pseudo-theta series and pseudo-Eisenstein series are series whose expressions closely resemble those of theta series and Eisenstein series. A typical example is given by the derivatives of Eisenstein series. In this paper, we will adapt the corresponding theory on quadratic spaces in \cite[Section 2]{Yuan1} to the setting of Hermitian spaces to align with our setting. We sketch the idea here.

Suppose $F$ is a totally real field and $E/F$ is a CM extension. Let $V$ be a Hermitian space over $E$, and let $\Phi=\Phi_S\otimes\Phi^S\in\mathcal{S}(V(\BA_E))$. The Weil representation of $\UU(\BA)$ gives the classical theta series
\begin{equation*}
    \theta(g,\Phi)=\sum_{x\in V(E)}r(g)\Phi(x),\quad g\in\UU(\BA).
\end{equation*}
Let $S$ be a finite set of non-archimedean places of $F$. Writing
\begin{equation*}
    r(g)\Phi(x)=r(g_S)\Phi_S(x)\,r(g^S)\Phi^S(x),
\end{equation*}
we obtain a pseudo-theta series by replacing the first local factor with a locally constant function $\Phi'_S(g_S,x)$:
\begin{equation*}
    A^{(S)}_{\Phi'}(g)=\sum_{x\in V(E)}\Phi'_S(g_S,x)\,r(g^S)\Phi^S(x).
\end{equation*}
We call it non-singular if $x\mapsto\Phi'_S(1,x)$ extends to a Schwartz function on $V(E_S)$. In that case its associated theta series is
\begin{equation*}
    \theta_{A^{(S)}}(g)=\sum_{x\in V(E)}r(g_S)\Phi'_S(1,x)\,r(g^S)\Phi^S(x),
\end{equation*}
and $A^{(S)}_{\Phi'}(g)=\theta_{A^{(S)}}(g)$ whenever $g_S=1$.

Now we start with the Siegel--Eisenstein series
\begin{equation*}
    E(s,g,\Phi)=\sum_{\gamma\in P(F)\backslash\UU(F)}\delta(\gamma g)^s r(\gamma g)\Phi(0),\quad g\in \UU(\BA),
\end{equation*}
where $\delta$ is the modulus function defined in \eqref{moduli character}. The non-constant part of $E(s,g,\Phi)$ has a Fourier expansion
\begin{equation*}
    E_*(s,g,\Phi)=\sum_{a\in F^\times}W_a(s,g,\Phi),
\end{equation*}
where $W_a(s,g,\Phi)$ is the Whittaker function. Again, let $S$ be a finite set of non-archimedean places of $F$. In
\begin{equation*}
    W_a(0,g,\Phi)=W_{a,S}(0,g,\Phi_S)W_a^S(0,g,\Phi^S),
\end{equation*}
if we replace $W_{a,S}(0,g,\Phi_S)$ by a locally constant function $B_{a,S}(g)$ of $(a,g)\in F^\times_S\times\UU(F_S)$, then we obtain a \textit{pseudo-Eisenstein series}
\begin{equation*}
    B_\Phi^{(S)}(g)=\sum_{a\in F^\times}B_{a,S}(g)W^S_a(0,g,\Phi^S),\quad g\in\UU(\BA).
\end{equation*}

Pseudo-Eisenstein series arise naturally in derivatives of Eisenstein series. In fact, for every non-archimedean $v$, the ``$v$-part"
\begin{equation*}
    \sum_{a\in F^\times}W'_{a,v}(0,g,\Phi)W^v_a(0,g,\Phi^v)
\end{equation*}
is a pseudo-Eisenstein series.

We say that a pseudo-Eisenstein series $B_\Phi^{(S)}(g)$ is \textit{non-singular} if for every $v\in S$, there exist $\Phi_v^+\in \mathcal{S}(V_v^+)$ and $\Phi_v^-\in \mathcal{S}(V_v^-)$ such that
\begin{equation*}
    B_{a,v}(1)=W_{a,v}(0,1,\Phi^+_v)+W_{a,v}(0,1,\Phi_v^-),\ \forall a\in F_v^\times.
\end{equation*}
Here $(V^+_v,V^-_v)$ is the set of (one or two) Hermitian spaces over $E_v$ with the same dimension and the same determinant as $V_v$. In this case, we form a linear combination of true Eisenstein series
\begin{equation*}
    E_B(s,g)=\sum_{\epsilon:S\rightarrow\{\pm\}}E(s,g,\Phi^\epsilon_S\otimes\Phi^S),\quad g\in\UU(\BA).
\end{equation*}
It approximates the original series in the sense that $B_\Phi^{(S)}(g)$ is equal to the non-constant part of $E_B(0,g)$ as long as $g_S=1$.

The key result in this ``pseudo'' theory is Lemma \ref{Key lemma}. Roughly speaking, this lemma states that if the sum of certain pseudo theta series and pseudo Eisenstein series is automorphic, then this sum is also equal to the sum of the corresponding theta series and Eisenstein series. In this way, we obtain a systematic method for controlling pseudo-theta series and pseudo-Eisenstein series.

The main value of this pseudo theory lies in its application to automorphic forms over general number fields. An important point is that classical modular forms over $\QQ$ possess many special properties that do not extend to general fields. A well-known example is the technique used by Richard Borcherds in \cite{Bor1,Bor2} proving the modularity of generating series on orthogonal Shimura varieties over $\QQ$. Borcherds has the following criterion of modularity.
\begin{theorem}\cite{Bor2}\label{Borcherds}
    A power series $f(\tau)=\sum_{n=0}^\infty a_nq^n$ converges to a holomorphic modular form of weight $k$ for $\SL_2(\ZZ)$ if and only if
    \begin{equation*}
        \sum_{n=0}^N a_nb_{-n}=0
    \end{equation*}
    for any modular form $g(\tau)=\sum_{m=-N}^\infty b_m q^m$ of weight $2-k$ for $\SL_2(\ZZ)$ that is holomorphic in $\mathcal{H}$ and meromorphic at infinity.
\end{theorem}
Using this criterion, Borcherds introduced a singular theta lifting, which is now widely known as the Borcherds lifting. This result has become a crucial technical tool in the Kudla program, playing an indispensable role in a series of works, c.f. \cite{BBK,BH,BHK+,HP,KRY,Kud3} for an incomplete list.

The particular obstruction argument above, which is phrased solely in terms of prescribed principal parts at cusps, does not extend directly to Hilbert modular forms of degree greater than one: by the Koecher principle, a Hilbert modular form cannot acquire poles supported only at the cusps. This does not assert that meromorphic Hilbert modular forms do not exist; rather, the classical principal-part criterion is no longer available in this form. The pseudo-series method provides a replacement suited to general totally real fields. It is already implicit in \cite{YZZ2}, was formulated more explicitly in \cite{YZ1}, and was developed systematically in \cite{Yuan1}.

In this series of papers, this pseudo theory also plays a crucial role. The final comparison in \cite{Guo2} largely depends on Lemma \ref{Key lemma}. In fact, compared with related work such as \cite{BH}, this pseudo theory compensates for the absence, over general number fields, of Borcherds lifting as well as the results proved in the appendix of \cite{BH}, thereby allowing us to carry out analogous arguments.

\subsection{Main theorem of the derivative series}
As we mentioned earlier, the \textit{derivative series} represents the ``analytic side" of our discussion. In one word, it is the holomorphic projection of the derivative of a mixed theta-Eisenstein series. This part mainly involves the transformation and computation of both theta series and Eisenstein series, as well as their derivatives, which is contained in Section \ref{Derivative series}. Here is a brief introduction.

Let $\BV$ be a totally positive definite incoherent Hermitian space over $\BA_E$. Fix an incoherent Hermitian subspace $\BW\subset \BV$ of dimension 1, and denote by $W^\perp(\BA_E)$ the orthogonal complement. Note that $W^\perp(\BA_E)$ is coherent, and we also denote by $W^\perp$ a Hermitian space over $E$.

For each Schwartz function $\Phi\in\mathcal{S}(\BV)$ invariant under an open compact subgroup $U$ of $\U(\BV_f)$, we have a \textit{mixed theta-Eisenstein series}
\begin{equation*}
    I(s,g,\Phi)=\sum_{\gamma\in P(F)\backslash \UU(F)}\delta(\gamma g)^s\sum_{x\in W^\perp}r(\gamma g)\Phi(x),
\end{equation*}
where $P$ is the subgroup of upper triangular matrices. The reason for using the term ``mixed'' is that for certain Schwartz functions, i.e., those Schwartz functions that can be written as $\Phi=\phi_1\otimes\phi_2$ with $\phi_1\in\mathcal{S}(W^\perp),\phi_2\in\mathcal{S}(\BW)$, this series happens to be a product of a theta series and an Eisenstein series.

The derivative $I'(0,g,\Phi)$ is an automorphic form of moderate growth, but it is neither holomorphic nor cuspidal. Since it need not be square-integrable, its cuspidal holomorphic projection is understood through Petersson pairings with holomorphic cusp forms (rather than as a literal Hilbert-space projection on the full space of automorphic forms). We denote this projection by $\Pr I'(0,g,\Phi)$ and call it the derivative series. A quasi-holomorphic projection $\Pr'I'(0,g,\Phi)$, defined coefficient by coefficient, will be used to calculate it explicitly.

Finally, we state the main theorem of this paper.
\begin{theorem}\label{main theorem of derivative series}
Assume the compatible choices of splitting characters, lattices, Haar measures, and local Schwartz functions specified in Sections \ref{Notation} and \ref{Choice of Schwartz function}. Then the derivative series $\Pr I'(0,g,\Phi)$ is an explicit finite linear combination of pseudo-theta series and pseudo-Eisenstein series. For the lattice-standard Schwartz function fixed in Section \ref{Choice of Schwartz function}, every local term occurring in this decomposition is computed explicitly in Section \ref{Explicit derivative}.
\end{theorem}

Since the explicit expression of the derivative series is rather lengthy, we do not write it out here. The proof of the first part of the main theorem follows from Theorem \ref{Projection of I'} and Theorem \ref{Projection of J'}. As we shall see, we use the quasi-holomorphic projection mentioned earlier to decompose the holomorphic projection into two parts and then derive explicit formulas for each of them separately.

The second part of the main theorem comes from the explicit computations carried out throughout Chapter \ref{Explicit derivative}, namely Lemma \ref{Explicit archimedean k series}, \ref{Explicit nonarchimedean k series}, \ref{Explicit archimedean}, \ref{Explicit nonarchimedean} and \ref{Explicit good intertwining}. It should be emphasized that these computations are all local in nature. Combined with the preceding discussion on pseudo-theta series and pseudo-Eisenstein series, this equivalently yields an explicit global formula.

One particularly important point to emphasize is that the logarithmic derivative of the L-function appears precisely in the explicit expression of the derivative series at $s=0$. More specifically, the local logarithmic derivative $L'(n+1,\eta_v^{n+1})/L(n+1,\eta_v^{n+1})$ appearing in Lemma \ref{Explicit nonarchimedean} is the $v$-component of one of the principal global terms in the modular-height formula. In order to obtain explicit results for this part, we need to extend many of the computational methods in \cite{Yuan1} to the higher-dimensional setting, and even introduce some new computational techniques. This constitutes one of the main innovations of our paper.

We now discuss some related works connected with our main theorem, as well as its applications. In fact, both the global study of derivatives of Eisenstein series and their holomorphic projections, and the local computation of explicit formulas for Whittaker functions, are of great importance in many areas, especially in work concerning special value formulas. One of the central problems in the Kudla program \cite{Kud4} is the so-called arithmetic Siegel–Weil formula, which (conjecturally) relate the first derivatives of Siegel Eisenstein series for unitary (resp. symplectic) groups with ``arithmetic theta series” formed from special cycles on integral models of unitary (resp. orthogonal) Shimura varieties. In several related works in recent years, such as \cite{LZ,LL1,LL2,Chen}, these topics arise simultaneously from both the local and global perspectives. For example, in \cite[Section 9]{LZ}, the local density is essentially a type of local Whittaker function. Our computation corresponds to one of the basic cases, i.e., when the fundamental matrix $T$ is a number (or equivalently when $n=1$ in their notations), and more general cases can then be obtained through the recursive relations on $n$ developed there.

Another more direct application of our results concerns another important problem in the Kudla program, namely the modularity of arithmetic generating series. Roughly speaking, this problem seeks to construct generating series on integral models taking values in arithmetic Chow groups. Related works include \cite{Kud3,KRY,BBK,BHK+,HP,Qiu} and others. This line of work furthermore has many concrete applications. For example, the proof of the arithmetic fundamental lemma in \cite{Zh3} relies on certain explicit computations carried out in \cite{BHK+}. In \cite{Qiu}, the author introduced a new approach to the modularity of arithmetic generating series, which involves extensive explicit computations of Eisenstein series and local Whittaker functions over general number fields. In particular, the constant $c_3$ appearing in our Theorem \ref{Projection of J'} provides the explicit value of the constant $\mathfrak{b}$ (up to a factor of 2) that was not made explicit in \cite[Thm 1.1.1, Remark 3.3.10]{Qiu}.

Finally, regarding our choice of Schwartz function in the main theorem, c.f. Section \ref{Choice of Schwartz function}, we provide an explanation for readers familiar with \cite{YZZ2,YZ1}. At \textit{each} finite place, our Schwartz function is taken to be the standard one, namely the characteristic function of a lattice in the Hermitian space, without imposing any “degeneracy” condition analogous to that in \cite[Section 5.2]{YZZ2}. Thus, compared with the computation in \cite{YZZ2,YZ1}, our expression is complicated since we cannot expect the disappearance of any local term.

\subsection{Organization of this series of papers}
In this paper, we define the derivative series which represents the ``analytic side" of our work, following the ideas we outlined above. The arrangement of this paper is as follows. In $\S$ \ref{Pseudo-theta series and Pseudo-Eisenstein series}, we give most of the definitions and terminologies used throughout this series of papers, and provide a rigorous theory of pseudo-theta series and pseudo-Eisenstein series. In $\S$ \ref{Derivative series}, we give a rigorous definition of the derivative series and complete the proof of the first half of the main Theorem \ref{main theorem of derivative series}. In $\S$ \ref{Explicit derivative}, we compute the explicit expressions of each pseudo-theta series and pseudo-Eisenstein series case by case at every local place, thereby completing the proof of the second half of the main theorem.

The second paper in this series is devoted to defining the height series which represents the ``arithmetic side" of our work. Roughly speaking, our height series is an automorphic form obtained from the arithmetic intersection numbers of the arithmetic generating series of divisors. In order to define it, we need to introduce the unitary Shimura variety and its integral model, the special divisors and their generating series, as well as their arithmetic counterparts. At the same time, we will also use the correspondence between quaternionic Shimura curves and unitary Shimura curves to derive the modular height formula for unitary Shimura curves, as well as the height formula for a CM point on them.

In the third and last paper of this series, we will carry out all the necessary computations and then compare the two series arising from the analytic side and the arithmetic side, thereby obtaining the modular height formula. We will find that most terms in these two series coincide, so in essence we also obtain a version of the arithmetic Siegel–Weil formula.

\subsection*{Acknowledgement}
The author is deeply grateful for the valuable assistance and meticulous guidance provided by professor Xinyi Yuan. Indeed, it is thanks to his previous work with Shou-Wu Zhang and Wei Zhang that the author has had the privilege to build upon it and make further contributions. He would like to thank his friend Weixiao Lu for much helpful advice. He thanks Ryan Chen, Yinchong Song, Liang Xiao and Roy Zhao for helpful communication. He also thanks Yifeng Liu and Wei Zhang for some useful suggestions. Finally, the author is grateful to the anonymous referee for so many valuable comments or suggestions to revise this paper.

\section{Pseudo-theta series and Pseudo-Eisenstein series}\label{Pseudo-theta series and Pseudo-Eisenstein series}
The goal of this section is to give a brief introduction to the pseudo-theta series and pseudo-Eisenstein series, and to introduce a very important lemma about these series, which we call \textit{the key lemma}. The original definition can be found in \cite[Section 2]{Yuan1}.

In the first half of this section, we will introduce various notations and terminologies used in this paper, especially focusing on those related to Hermitian spaces and Weil representations. Later, we recall the notion of Whittaker functions, theta series and Eisenstein series in our setting, then introduce the general theory of pseudo-theta series and pseudo-Eisenstein series.

\subsection{Hermitian space and Weil representation}\label{Notation}
We list some notations and terminologies we will use in this series of papers. Almost all the notations are compatible with the notations in \cite{YZZ2,YZ1,Yuan1}. Throughout this paper, $F$ is a totally real field. $E/F$ is a quadratic totally imaginary extension, i.e., $E$ is a CM field. Denote by $\BA_F$ the ring of adeles of $F$, which we sometimes use $\BA$ for convenience, and $\Af$ the ring of finite adeles of $F$. We use $v$ for a place of $F$. For convenience, we use $\BA_F^v$ for the restricted product of all local fields of $F$ omitting the $F_v$-term.

Similarly, we can define these notations for $E$. In this series of papers, if there is no subscript of the absolute value $|\cdot|$, we mean the absolute value of $F_v$. Otherwise, we will use $|\cdot|_{E_v}$ for the absolute value of $E_v$. By class field theory, note that there is a natural quadratic character $\eta$ of $\BA^\times$ corresponding to $E/F$.

Throughout the paper, $L(s,\eta)=\prod_vL(s,\eta_v)$ denotes the completed Hecke $L$-function with the local normalizations fixed below, and $L_f(s,\eta)$ denotes its finite part.

We also use the notation $\avint_G$ for the normalized integral $\frac{1}{\vol(G(F_v))}\int_{G(F_v)}$, here $G$ is any algebraic group over $F$. This notation is also valid when we consider the integral over $G(\BA)$.

While for Hermitian spaces, throughout this paper, we always denote the Hermitian pairing by $\langle\cdot,\cdot\rangle$ and the Hermitian norm by $q(\cdot)$.

\subsubsection*{Additive character and measure}

For a place $v$ of the totally real field $F$, write $\Fn=F_v$. Thus an archimedean $F_v$ is always $\RR$. We normalize
\begin{equation*}
 \psi_v(x)=e^{2\pi ix}\quad(F_v=\RR),
 \qquad
 \psi_v=\psi_{\QQ_p}\circ\tr_{F_v/\QQ_p}\quad(v\nmid\infty),
\end{equation*}
where $\psi_{\QQ_p}(x)=e^{-2\pi i\iota(x)}$. The character used on $E_v$ is
$\psi_{E_v}=\psi_v\circ\tr_{E_v/F_v}$; in particular, if $E_v=\CC$, then
$\psi_{E_v}(z)=e^{4\pi i\Re(z)}$.

We normalize the measure $dx$ on $\Fn$ to be the unique Haar measure on $\Fn$ self-dual
with respect to $\psi$ in the sense that the Fourier transform
\begin{equation*}
    \widehat{\Phi}(y):=\int_{\Fn}\Phi(x)\psi(xy)dx
\end{equation*}
satisfies the inversion formula $\widehat{\widehat{\Phi}}(x)=\Phi(-x)$. The measures are
determined explicitly as follows:

If $\Fn=\RR$, then $dx$ is the usual Lebesgue measure.

For later use on $E_v=\CC$, the self-dual measure for $\psi_{E_v}(z)=e^{4\pi i\Re z}$ is twice the usual Lebesgue measure.

If $\Fn$ is non-archimedean, then $\vol(\mathcal{O}_{\Fn})=|d_{\Fn}|^{\frac{1}{2}}$. Here
$\mathcal{O}_{\Fn}$ is the ring of integers and $d_{\Fn}$ is the different of $\Fn$ over $\QQ_p$.

We also normalize the Haar measure $d^\times x$ on $\Fn^\times$ by
\begin{equation*}
    d^\times x=\zeta_{\Fn}(1)|x|^{-1}dx.
\end{equation*}
Here $\zeta_{\Fn}(s)=(1-N^{-s})^{-1}$ if $\Fn$ is non-archimedean whose residue field has
$N$-elements, $\zeta_{\RR}(s)=\pi^{-s/2}\Gamma(s/2)$, and $\zeta_{\CC}(s)=2(2\pi)^{-s}\Gamma(s)$. This normalization implies if $\Fn$ is non-archimedean, then $\vol(\mathcal{O}_{\Fn}^\times,d^\times x)=\vol(\mathcal{O}_{\Fn},dx)$.

\subsubsection*{Coherent and incoherent Hermitian space}
In this article, we only discuss non-degenerate Hermitian space. We refer to \cite{Ser} and \cite{Gro2} for general properties of Hermitian spaces. Let $F$ be a totally real field, and let $E/F$ be a totally imaginary quadratic extension. Denote by $\tau$ the non-trivial involution of $E$ over $F$. For a Hermitian space $(V,q)$ over $E$, it determines a unique element $d_V\in F^\times/\Nm_{E/F}(E^\times)$ which is called the \textit{Hermitian determinant}.

Note that similarly, we can also define Hermitian spaces when $E/F$ is the \'{e}tale quadratic algebra $F\oplus F$ with involution $\tau(a,b)=(b,a)$. In this case, $d_V=1$ obviously. Indeed, in this case the Hermitian space is a free $E$-module of the form $V=W\oplus W'$, where $W$ is a vector space over $F$ and $W'$ is its dual space. Then if we set $v=w+u'$ with $w\in W$ and $u'\in W'$, the Hermitian pairing is defined by $q(v)=\langle w,u'\rangle\in F$.

For convenience, in the further discussion, for any $t\in F^\times$, we denote by $V_t\subset V$ the subset of those elements $v\in V$ with $q(v)=t$. This notation will be used throughout this series of papers.

It is clear that all the notations are valid when passing to local case.
For convenience, for any place $v$ of $F$, denote by $\epsilon(V_v)=(d_{V,v},D_{E/F})_v$ the local Hasse invariant of $(V_v,q)$, where $d_{V,v}$ is the local Hermitian determinant and $D_{E/F}$ is the relative discriminant.
Note that if $v$ is split in $E/F$, $(V_v,q)$ is unique with local Hasse invariant 1.
If $v$ is a non-archimedean place which is nonsplit in $E/F$, $E_v$ is a quadratic field over $F_v$ and the group $F_v^\times/\Nm(E_v^\times)$ has order two. Thus, in
this case, there are exactly two local Hermitian spaces with opposite Hasse invariant.
In many references, the case where $\epsilon(V_v)=1$ is referred to as "split," while the case where $\epsilon(V_v)=-1$ is referred to as "nonsplit."
When $F_v\cong\RR$ and $E_v\cong \CC$, there is a further invariant of a Hermitian space called signature. If the signature is $(r,s)$, then $d_{V,v}=(-1)^s$.

There is a well know local-to-global principle for Hermitian spaces over CM field.
Fix a dimension $n\ge1$, denote by $\BV$ an adelic Hermitian space over $\BA_E$, and
assume that the local Hasse invariant $\epsilon(\BV_v)=1$ for almost all $v$. Then a
necessary and sufficient condition for the existence of a global Hermitian space
$(V_0,q)$ such that $(V_0(\BA_E),q)\cong(\BV,q)$ is
\begin{equation*}
    \prod_v \epsilon(\BV_v)=1,
\end{equation*}
and if this relation holds, $(V_0,q)$ is unique up to isomorphism. Set $\epsilon(\BV)=\prod_v\epsilon(\BV_v)\in\{\pm1\}$. We call $\BV$ coherent if $\epsilon(\BV)=1$ and incoherent if $\epsilon(\BV)=-1$.

In our setting, we fix an incoherent Hermitian space $(\BV,q)$ of dimension $n+1$ over $\BA_E$, which is positive definite at each archimedean place, and we also fix a real place $\iota$ of $F$. Then there is a unique coherent Hermitian space $(V,q)$ over $E$
such that $V_v\cong\BV_v$ for any place $v\ne\iota$, while $V_\iota$ has signature $(n,1)$. We call such $V$ the \textit{nearby coherent Hermitian space} of $\BV$ with respect to $\iota$. Similarly, we can talk about the nearby coherent Hermitian space $V(v)$ with
respect to any place $v$ which is nonsplit in $E/F$.

\subsubsection*{Decomposition of Hermitian space}
Following the notations above, let $V/E$ be a Hermitian space that has signature $(n,1)$ at $\iota$ and signature $(n+1,0)$ at all other archimedean places. Denote by $\BV$ the incoherent totally positive definite Hermitian space associate to $V$. We further choose an orthogonal decomposition
\begin{equation}\label{standard decomposition of incoherent Hermitian space}
    \BV=\BW\oplus W^\perp(\BA_E),
\end{equation}
where $\BW$ is an incoherent totally positive definite Hermitian subspace of dimension 1. Thus, the orthogonal complement is coherent. The choice of $\BW$ is given as follows. Using the local-to-global principle above, there is a unique incoherent
Hermitian space $\BW$ of dimension one up to isomorphism, such that $d_{\BW,v}=d_{\BV,v}$ for all places.
Let $W^\perp$ be the coherent orthogonal complement, then we have $\epsilon(W^\perp_v)=1$ at all places. The simplicity of $W^\perp$ is useful for the later computation.

Similarly, there is also an orthogonal decomposition
\begin{equation*}
    V=W\oplus W^\perp
\end{equation*}
for the coherent Hermitian space. Here $W$ is also the nearby coherent Hermitian space of $\BW$ with respect to the archimedean place $\iota$.
We use $q$ for the Hermitian form of $\BV$ or $V$, and for convenience, $q$ also means the Hermitian form restricted to $\BW$ or $W^\perp$.

In the subsequent discussion we fix an incoherent totally positive definite Hermitian subspace $\BV_1\subset\BV$ of dimension $2$ containing $\BW$ and satisfying $d_{\BV_1}=d_{\BW}=d_{\BV}$. Equivalently, we choose a coherent Hermitian line in $W^\perp$ with the required determinant; the existence of this compatible auxiliary choice is part of the standing data. Denote by $V_1$ the nearby Hermitian space of $\BV_1$ with respect to the archimedean place $\iota$. There are orthogonal decompositions
\begin{equation}\label{a further decomposition}
    \BV=\BV_1\oplus V_1^\perp(\BA_E),\quad \BV_1=\BW_1\oplus W_1^\perp(\BA_E),
\end{equation}
where $V_1^\perp$ is the coherent orthogonal complement of $\BV_1$ in $\BV$ and $W_1^\perp$ is the coherent orthogonal complement of $\BW_1=\BW$ in $\BV_1$. We can also define the decomposition of nearby coherent Hermitian spaces with respect to $\iota$.

\subsubsection*{Lattice in Hermitian space}
Keep all the notations above, and denote by $\Lambda^\vee_v=\{x\in\BV_v:\langle x,\Lambda_v\rangle\subset\mathcal{O}_{E_v}\}$ the dual lattice of $\Lambda_v$. We introduce the following definition for an $\mathcal{O}_{E_v}$ lattice $\Lambda_v\subset\BV_v$. See also \cite[Section 4.4]{Qiu}:
\begin{itemize}
    \item $\Lambda_v$ is \textit{self-dual} if $\Lambda_v^\vee=\Lambda_v$;
    \item $\Lambda_v$ is \textit{$\varpi_{E_v}$-modular} if $\Lambda_v^\vee=\varpi^{-1}_{E_v}\Lambda_v$;
    \item  $\Lambda_v$ is \textit{almost $\varpi_{E_v}$-modular} if $\Lambda_v^\vee\subset\varpi^{-1}_{E_v}\Lambda_v$ and the inclusion is of colength 1.
\end{itemize}
Here $\varpi_{E_v}$ is the uniformizing parameter of $E_v$.

\begin{remark}\label{dual by trace remark}
    It is important to emphasize that at ramified places, our notation is the same as the notation in \cite{RSZ1}, \cite{BH} and \cite{Qiu}, but is not the same as the one in \cite{LL2}. We remark that the $\varpi_{E_v}$-modular lattice $\Lambda_v$ is the one that is more compatible with Fourier transforms. Indeed, it is not hard to check that
    \begin{equation*}
        \varpi_{E_v}^{-1}\Lambda_v=\{x\in \BV_v\Big|\tr_{E_v/F_v}(\langle x,y\rangle)\in\mathcal{O}_{F_v},\ \forall y\in\varpi_{E_v}^{-1}\Lambda_v\}.
    \end{equation*}
    This is the fundamental reason why we utilize this particular class of Hermitian lattices in the further discussion. Meanwhile, in some other literature, such as \cite{LL2}, the author uses $\tr_{E_v/F_v}(\langle \cdot,\cdot\rangle)$ to define the dual of a Hermitian lattice. Then $\varpi_{E_v}^{-1}\Lambda_v$ becomes their "self-dual lattice" (after taking the trace map).
\end{remark}

For any incoherent Hermitian space $\BV$ as above, we will consider a Hermitian lattice as follow.
\begin{definition}\label{Choice of lattice}
Define a Hermitian lattice
\begin{equation*}
    \Lambda=\prod_{v\nmid\infty}\Lambda_v\subset\BV_f
\end{equation*}
such that $\Lambda_v$ is self-dual if $v$ is unramified in $E/F$ and is $\varpi_{E_v}$-modular (resp. almost $\varpi_{E_v}$-modular) if $v$ is ramified in $E/F$ with $2\nmid n$ (resp. $2|n$).
\end{definition}

\begin{lemma}\label{Existence of lattice}
Assume that $\epsilon(\BV_v,q)=1$ at all inert places.
\begin{enumerate}
 \item If $n$ is odd, then a lattice as in Definition \ref{Choice of lattice} exists if and only if $(n+1)/2$ and $[F:\QQ]$ are both odd.
 \item If $n$ is even, then a lattice as in Definition \ref{Choice of lattice} exists.
\end{enumerate}
\end{lemma}
\begin{proof}
This follows from the local classification of the ramified lattices and the product formula for the discriminants; see \cite[Remark 5.1.2]{Qiu}.
\end{proof}

The statement concerns the full adelic lattice with all the local conditions in Definition \ref{Choice of lattice}, not merely the existence of an almost $\varpi_{E_v}$-modular lattice at one place.

Moreover, recall the orthogonal decompositions \ref{standard decomposition of incoherent Hermitian space} and \ref{a further decomposition}. At every finite place $v$, choose a representative $d_{\BV,v}\in F_v^\times$ of the determinant class and a vector $e_v\in\BW_v$ with
\begin{equation*}
 q(e_v)=d_{\BV,v}\in\mathcal O_{F_v}^\times.
\end{equation*}
Set $\Lambda_{\BW_v}=\mathcal O_{E_v}e_v$. At every unramified place we assume a compatible orthogonal decomposition
\begin{equation}\label{decomposition of inert Hermitian lattice}
 \Lambda_v=\Lambda_v^\perp\oplus\Lambda_{\BW_v}.
\end{equation}
At a ramified place, the local classification in \cite[Lemma 5.1.1 and Section 5.1.3]{Qiu} allows us to choose $\Lambda_{1,v}=\Lambda_v\cap\BV_{1,v}$ as a rank-two $\varpi_{E_v}$-modular summand, with $\Lambda_v=\Lambda_{1,v}\oplus\Lambda_{1,v}^\perp$ and $\BW_v\subset E_v\Lambda_{1,v}$. Choosing $e'_v\in\BV_{1,v}$ orthogonal to $e_v$ with unit norm, we then have
\begin{equation}\label{decomposition of ramified Hermitian lattice}
 \varpi_{E_v}(\mathcal O_{E_v}e_v\oplus\mathcal O_{E_v}e'_v)
 \subset\Lambda_{1,v}\subset
 \mathcal O_{E_v}e_v\oplus\mathcal O_{E_v}e'_v,
\end{equation}
with both inclusions of colength one. We put
\begin{equation*}
 \Lambda_v^\perp=\mathcal O_{E_v}e'_v\oplus\Lambda_{1,v}^\perp.
\end{equation*}
In particular, $\Lambda_v\neq\Lambda_v^\perp\oplus\Lambda_{\BW_v}$ at a ramified place.

Note that using the result in \cite[Appendix 1.B]{GKT}, we can also compute some local Weil index using the local lattice. A direct result is that if $v$ is split in $E/F$, then $\gamma(\BV_v,q)=1$ since there is a self-dual lattice; if $v$ is inert in $E/F$, then $\gamma(\BV_v,q)=\epsilon(\BV_v,q)$; if $v$ is ramified in $E/F$, then we also have $\gamma(\BV_v,q)=1$ for the same reason.

\subsubsection*{Volume of local quadratic extension}
Keep all the notations above, we now normalize the Haar measures on the quadratic extension $E_v/F_v$ based on the normalization of $F_v$ above.

To see the Haar measure, let $E_v$ be a quadratic \'{e}tale algebra extension over $\Fn$ with $q:E_v\longrightarrow\Fn$
the reduced norm. One has one of the following:
\begin{itemize}
    \item $E_v=\Fn\oplus\Fn$ and $q(x_1,x_2)=x_1x_2$, in this case we say $E_v$ is a split quadratic extension;
    \item $E_v$ is a quadratic field extension over $\Fn$ and $q=\Nm_{E_v/\Fn}$, in this case
    we say $E_v$ is a nonsplit quadratic extension.
\end{itemize}
In both case, we endow $E_v$ with the self-dual Haar measure $dx$ with respect to the local quadratic space $(E_v,q)$, and endow $E_v^\times$ with a Haar measure $d^\times x$ with
\begin{equation*}
    d^\times x=\zeta_{E_v}(1)|q(x)|^{-1}dx.
\end{equation*}
Here $\zeta_{E_v}(s)=\zeta_{\Fn}(s)^2$ in the split case, while $\zeta_{E_v}(s)$ is the zeta
function of $E_v$ in the nonsplit case.

We also endow the subgroup
\begin{equation*}
    E_v^1:=\{h\in E_v|q(h)=1\}
\end{equation*}
with the Haar measure $dh$ determined by the exact sequence
\begin{equation*}
    1\longrightarrow E_v^1\longrightarrow E_v^\times\longrightarrow q(E_v^\times)\longrightarrow 1.
\end{equation*}

Following the result in \cite[Section 1.6.2]{YZZ2}, we also have the following explicit results. In the split case, the Haar measures on $E_v$ and $E_v^\times$ are compatible with
the Haar measure $dx$ and $d^\times x$ respectively. In the nonsplit case, if $\Fn$ is
non-archimedean, still denote by $d=d_{\Fn}$ the different of $\Fn$ over $\QQ_p$, and by
$D=D_{E_v/\Fn}\in\mathcal{O}_{\Fn}$ the discriminant of $E_v/\Fn$, we have
\begin{equation*}
    \vol(\mathcal{O}_{E_v},dx)=\vol(\mathcal{O}_{E_v}^\times,d^\times x)=|D|^{\frac{1}{2}}|d|.
\end{equation*}
Furthermore,
\begin{equation*}
    \vol(E_v^1)=\left\{
    \begin{aligned}
        \nonumber
        &2\ \ \ \mathrm{if}\ \Fn=\RR\ \mathrm{and}\ E_v=\CC,\\
        &|d|^\frac{1}{2}\ \ \ \mathrm{if}\ E_v/\Fn\ \mathrm{is\ nonsplit\ and\ unramified}.\\
        &2|D|^{\frac{1}{2}}|d|^\frac{1}{2}\ \ \ \mathrm{if}\ E_v/\Fn\ \mathrm{is\ ramified}.\\
    \end{aligned}
    \right.
\end{equation*}
With these local normalization, we also have the corresponding normalization in global case by taking product. Especially, we have
\begin{equation*}
    \vol(E^1\backslash E^1(\BA))=2L(1,\eta).
\end{equation*}

Fix the self-dual Haar measure on $\BV_v$, by our choice of the Hermitian lattice $\Lambda_v$ in Definition \ref{Choice of lattice}, we conclude that for any unramified prime,
\begin{equation*}
    \vol(\Lambda_v)=|d_v|^{n+1},\ \vol(\Lambda_{W_v})=|d_v|,
\end{equation*}
For the explicit formulas below, we impose the standing hypothesis $|d_v|=1$ at every place ramified in $E/F$ (as restated in Section \ref{Choice of Schwartz function}). Under this hypothesis, at a ramified place we have
\begin{equation*}
    \vol(\Lambda_v)=\left\{
    \begin{aligned}
        \nonumber
        &N_v^{-(n+1)}\ \ \ 2\nmid n,\\
        &N_v^{-(n+\frac{1}{2})}\ \ \ 2|n,\\
    \end{aligned}
    \right.\ \vol(\Lambda_{W_v})=N_v^{-\frac{1}{2}}.
\end{equation*}
These results can be checked easily from definition. We will use this numerical result in the later explicit computation.

\subsubsection*{Schwartz functions}
We give a brief introduction of those Schwartz functions we will use in this paper. Here we again denote by $(V,q)$ a Hermitian space over $E$.

For non-archimedean places, in general, we make no additional assumption, i.e., we consider those locally constant and compactly supported functions. Nonetheless, for the explicit computation in the later discussion, we will make some extra restrictive assumptions for convenience. See Section \ref{Choice of Schwartz function}.

Let $v$ be an archimedean place of $F$, then $F_v=\RR$ and $E_v=\CC$. The standard Schwartz function $\Phi_v\in\mathcal{S}(V(E_v))$ is the Gaussian function
\begin{equation*}
    \Phi_v(x)=e^{-2\pi q(x)}.
\end{equation*}
For all explicit formulas in this series we impose the standing condition that the archimedean components of $\Phi$ are the standard Gaussians. This is a restriction on the chosen Schwartz datum, made because it is the datum compatible with the holomorphic generating series used in the later papers.

\subsubsection*{Subgroup of $\U(1,1)$}
Let $\U(1,1)$ be the unitary group over $F$ of the two-dimensional standard skew-Hermitian space over $E$. For convenience, we introduce the following matrix notation:
\begin{equation}\label{moduli character}
    m(a)=\matrixx{a}{}{}{\bar{a}^{-1}},\ n(b)=\matrixx{1}{b}{}{1},\ w=\matrixx{}{1}{-1}{}.
\end{equation}
Here $a\in E^\times$ and $b\in F$, and $\bar{a}$ means $\tau(a)$. We also denote by $P$ the subgroup of upper triangular matrices, $M$ the diagonal subgroup generated by $m(a)$, and $N$ the standard unipotent subgroup generated by $n(b)$. Clearly, the isomorphism $N\cong \BG_{a,F}$ induces an additive character and Haar measure on $N(\BA_F)$.

For any place $v$ of $F$, the character
\begin{equation*}
    \delta_v:P(\Fn)\longrightarrow\RR^\times,\ \matrixx{a}{b}{}{\bar{a}^{-1}}\mapsto\left|a\right|^{\frac{1}{2}}_{E_v}=\left|\Nm_{E_v/F_v}(a)\right|^{\frac{1}{2}}_{F_v}
\end{equation*}
extends to a function $\delta_v:\UU(\Fn)\longrightarrow\RR^\times$ by Iwasawa
decomposition
\begin{equation*}
    \UU(F_v)=N(F_v)M(F_v)K(F_v).
\end{equation*}
Here $K(F_v)$ is the intersection of $\UU(F_v)$ with the standard maximal compact subgroup of $\GL_2(E_v)$. We immediately know that for an archimedean place $v$, $K(F_v)$ is the subgroup of form
\begin{equation*}
    [k_1,k_2]:=\frac{1}{2}\matrixx{k_1+k_2}{-ik_1+ik_2}{ik_1-ik_2}{k_1+k_2},
\end{equation*}
where $k_1,k_2\in E_v\cong\CC$ are of norm 1.

For example, we have
\begin{equation*}
    \delta_v(wn(b))=\left\{
         \begin{aligned}
         \nonumber
         & 1\quad b\in \mathcal{O}_{\Fn},\\
         & |b|^{-1}\quad \mathrm{otherwise}.
        \end{aligned}
        \right.
\end{equation*}
Indeed, if $v(b)<0$, we have
\begin{equation*}
    wn(b)=n(-b^{-1})m(b^{-1})wn(-b^{-1})w,
\end{equation*}
where the matrix $wn(-b^{-1})w$ lies in $K(F_v)$.

For the global field $F$, the product $\delta=\prod_v\delta_v$ gives a function on $\UU(\BA)$.

We should remind the reader that, our modulus character is compatible with the one in \cite{YZZ2,YZ1,Yuan1}, but is not compatible with the one in \cite{Qiu}.

\subsubsection*{Weil representation}
We give a brief introduction to the Weil representation associated with $\UU(F_v)$. For a local Hermitian space $V_v$, the group $\UU(F_v)\times\U(V_v)$ acts on the Schwartz space $\mathcal S(V_v)$, not on the underlying vector space itself. For any $\Phi\in\mathcal S(\BV_v)$, the action is given as follows:
\begin{itemize}
    \item $r(h)\Phi(x)=\Phi(h^{-1}x)$,\quad $h\in\U(\BV_v)$;
    \item $r(m(a))\Phi(x)=\chi_{\BV,v}(a)|a|_{E_v}^{\dim \BV/2}\Phi(ax)$,\quad $a\in E_v^\times$;
    \item $r(n(b))\Phi(x)=\psi(bq(x))\Phi(x)$,\quad $b\in F_v$;
    \item $r(w)\Phi=\gamma(\BV_v,q)\widehat{\Phi}$, \quad $w=\matrixx{}{1}{-1}{}$.
\end{itemize}
We now specify the splitting characters. As part of the standing Weil data, fix a unitary Hecke character $\chi:E^\times\backslash\BA_E^\times\to\CC^\times$ with $\chi|_{\BA_F^\times}=\eta$, and choose it unramified at every finite place at which $E/F$ and the lattice data are unramified. For a Hermitian space $V$ of dimension $d$, take $\chi_V=\chi^d$; equivalently, one may use any compatible family satisfying $\chi_V|_{\BA_F^\times}=\eta^d$ and $\chi_{V_1\oplus V_2}=\chi_{V_1}\chi_{V_2}$. Merely prescribing the restriction to $\BA_F^\times$ does not determine $\chi_V$ uniquely, while the unramified condition is part of the spherical datum used below. The constant $\gamma(\BV_v,q)$ is the Weil index. Finally,
$\widehat{\Phi}$ denotes the Fourier transform
\begin{equation*}
 \widehat{\Phi}(x)=\int_{\BV_v}\Phi(y)\,
 \psi_v\!\left(\tr_{E_v/F_v}\langle x,y\rangle\right)dy,
\end{equation*}
where $dy$ is the self-dual Haar measure on $\BV_v$ for this pairing.

In this article, we will use $r(g,h)$ for Weil representation on different Hermitian spaces, including $\BW_v$, $W^\perp_v$ and $\BV_v$ with $\BV=\BW\oplus W^\perp(\BA)$ from the fixed decomposition. We can clearly see from the choice of the Schwartz function which space the Weil representation we are using.

\subsection{Theta series and Eisenstein series}\label{Theta}
In this subsection, we define the theta series and Eisenstein series that we will use in this paper. Most of the definitions below apply in greater generality.

For a totally real field $F$ and $E/F$ a CM field, let $(V,q)$ be a Hermitian space over $E$. Recall that we assume $V$ has signature $(n,1)$ at a fixed archimedean place $\iota$ and has signature $(n+1,0)$ at all other archimedean places. Denote by $\BV$ the incoherent totally positive definite Hermitian space associate to $V$. Recall that we have a decomposition $\BV=\BW\oplus W^\perp(\BA_E)$, where $\BW$ is an incoherent totally positive definite Hermitian subspace of dimension 1, and denote by $W$ the
nearby coherent Hermitian space of $\BW$ with respect to $\iota$.

\subsubsection*{Whittaker function}
Recall the notion of Whittaker function. For any $a\in F$, the $a$-th Whittaker function
is defined by
\begin{equation*}
    W_a(s,g,\Phi)=\int_\BA\delta(wn(b)g)^s r(wn(b)g)\Phi(0)\psi(-ab)db.
\end{equation*}
We also have the local Whittaker integrals
\begin{equation*}
    W_{a,v}(s,g,\Phi_v)=\int_{F_v}\delta_v(wn(b)g)^s r(wn(b)g)\Phi_v(0)\psi(-ab)db.
\end{equation*}
Here $g\in\UU(F_v)$ and $\Phi_v\in\mathcal{S}(\BV_v)$ in the local integral; the global integral uses $\Phi\in\mathcal S(\BV)$. Note that we will also consider the Whittaker function for $\phi_2\in\mathcal{S}(\BW)$ later. These two functions have a meromorphic continuation to all $s\in\CC$.

We also have the following equations
\begin{equation}\label{action on Whittaker function}
    W_{a,v}(s,n(b)g,\Phi_v)=\psi(ab)W_{a,v}(s,g,\Phi_v),\ b\in \Fn,
\end{equation}
\begin{equation*}
    W_{a,v}(s,m(a')g,\Phi_v)=|a'|_{E_v}^{1-\frac{\dim\BV}{2}-\frac{s}{2}}\chi_{\BV,v}(a')W_{\Nm(a')a,v}(s,g,\Phi_v),\ a'\in E_v^\times.
\end{equation*}
Note that the first equation is obvious, while the second equation is deduced from the
relation $wn(b)m(a)=m(\bar{a}^{-1})wn(\Nm(a)^{-1}b)$.

\subsubsection*{Theta series}

We form a theta series associated with $W^\perp$. For any $\phi_1\in\mathcal{S}(W^\perp(\BA_E))$,
define a function on $\UU(F)\backslash\UU(\BA)\times \U(W^\perp)\backslash \U(W^\perp(\BA_E))$:
\begin{equation*}
    \theta(g,h,\phi_1)=\sum_{x\in W^\perp} r(g,h)\phi_1(x),\quad (g,h)\in\UU(\BA)\times \U(W^\perp(\BA_E)).
\end{equation*}
Note that when  we take $h=\id$, we use $r(g)$ and $\theta(g,\phi_1)$ for simplicity.
We also call
\begin{equation*}
    \theta_0(g,\phi_1)=r(g)\phi_1(0)
\end{equation*}
the constant term of theta series.

Analogous theta series are defined for every coherent global Hermitian space, in particular for $V$ and $W$. When the Schwartz data are initially given on the incoherent adelic spaces $\BV$ or $\BW$, a global theta series is formed only after replacing the relevant local component so that the resulting adelic Hermitian space is coherent. We use the notation $\theta(g,h,\Phi)$ for that coherent realization.

\subsubsection*{Eisenstein series}

Following the idea in \cite{YZZ2}, where the origin is in the work of Gross and Zagier \cite{GZ}, we will consider the incoherent Eisenstein series associated with $\BW$.
Let $\phi_2\in\mathcal{S}(\BW)$ be a Schwartz function. The associated Siegel--Eisenstein series is given by
\begin{equation*}
    E(s,g,\phi_2)=\sum_{\gamma\in P(F)\backslash\UU(F)}\delta(\gamma g)^s r(\gamma g)\phi_2(0),\quad g\in\UU(\BA).
\end{equation*}
It has a meromorphic continuation to $s\in\CC$. Moreover, it is holomorphic at $s=0$.

By the standard theory, we have
\begin{equation*}
    E(s,g,\phi_2)=\delta(g)^s r(g)\phi_2(0)+W_0(s,g,\phi_2)+\sum_{a\in F^\times}W_a(s,g,\phi_2).
\end{equation*}
Note that we also have the constant term of Eisenstein series
\begin{equation*}
    E_0(s,g,\phi_2)=\delta(g)^s r(g)\phi_2(0)+W_0(s,g,\phi_2).
\end{equation*}

Later, in the computation of \ref{J sequence}, we will also consider the incoherent Eisenstein series associated with $\BV$: for any Schwartz function $\Phi\in\mathcal{S}(\BV)$, the associated Eisenstein series is given by
\begin{equation*}
    E(s,g,\Phi)=\sum_{\gamma\in P(F)\backslash\UU(F)}\delta(\gamma g)^s r(\gamma g)\Phi(0),\quad g\in\UU(\BA).
\end{equation*}
It is not hard to distinguish Eisenstein series associated with different spaces by
checking its Schwartz function. In this article, we will always use $\Phi$ for a
Schwartz function in $\BV$ and $\phi_2$ for a Schwartz function in $\BW$.

For convenience, for any $t\in F$, denote by $I_t(s,g,\Phi)$ the $\psi_t$-Whittaker function of $I(s,g,\Phi)$. Later in Theorem \ref{Projection of I'} and Theorem \ref{Projection of J'}, we will present our result from a holistic perspective, following the notations in \cite{YZZ2,YZ1,Yuan1}. But we remind the reader that the same result holds for each $\psi_t$-component when $t\in F^+$. See \cite{Qiu} for example. In our final comparison between the derivative series and height series, we need the result of each $I_t(s,g,\Phi)$ for any $t\in F^+$, since our height series is a sum parametrized by such $t$.

\subsection{Pseudo-theta series and pseudo-Eisenstein series}\label{Pseudo theta and Eisenstein}
In this subsection, we define pseudo-theta series and pseudo-Eisenstein series. We mainly refer to \cite[Section 6]{YZ1} and \cite[Section 2]{Yuan1}, and note that we consider the group $\UU$ here, not $\GL_2$. Nonetheless, this makes no essential difference for the result and proof. One can compare the definition here with the usual theta series and Eisenstein series defined in Section \ref{Theta}.

\subsubsection*{Pseudo-theta series}
Fix a quadratic field extension $E/F$,
let $V$ be a Hermitian space over $E$, and $V_0\subset V_1\subset V$ be two subspaces over $E$ with induced Hermitian forms. Moreover,
we allow $V_0=\emptyset$. Denote by $S$ a finite set of non-archimedean places
of $F$, and $\Phi^S\in\mathcal{S}(V(\BA_E^S))$ be a Schwartz function with standard
infinite components.

A pseudo-theta series is a series of the form
\begin{equation*}
    A_{\Phi'}^{(S)}(g)=\sum_{x\in V_1(E)\setminus V_0(E)}\Phi'_S(g_S,x)r_V(g^S)\Phi^S(x),\quad g\in\UU(\BA).
\end{equation*}
We explain these notations as follows:
\begin{itemize}
    \item The Weil representation is $r_V$ attached to $V$.
    \item $\Phi'_S(g,x)=\prod_{v\in S}\Phi'_v(g_v,x_v)$ as local product, with
    \begin{equation*}
        \Phi'_v(\cdot,\cdot):\UU(F_v)\times(V_1-V_0)(E_v)\longrightarrow\CC
    \end{equation*}
    a smooth locally constant function for each $v\in S$. Here the smoothness
    means the right translation of some open compact subgroup $K_v\subset\UU(F_v)$ is invariant.
    \item For any $v\in S$ and $g\in\UU(F_v)$, the support of $\Phi'_v(g,\cdot)$ in $(V_1-V_0)(E_v)$ is bounded, which makes the sum
    convergent.
\end{itemize}

The pseudo-theta series $A^{(S)}$ is called non-degenerate (resp. degenerate) if
$V_1=V$ (resp. $V_1\ne V$). It is called non-singular if for each $v\in S$, the
local component $\Phi'_v(1,x)$ can be extended to a Schwartz function on $V_1(E_v)$. One remark is that in our case, all the pseudo-theta
series are non-singular. For convenience, we sometimes omit the subscript $\Phi'$
of the pseudo-theta series.

There are two usual theta series associated with a non-singular pseudo-theta series
$A^{(S)}_{\Phi'}$. The theta series
\begin{equation*}
 \theta_{A,1}(g)=\sum_{x\in V_1(E)}r_{V_1}(g_S)\Phi'_S(1,x)\,r_{V_1}(g^S)\Phi^S(x)
\end{equation*}
is called the outer theta series associated with $A^{(S)}_{\Phi'}$. Similarly, the theta series
\begin{equation*}
 \theta_{A,0}(g)=\sum_{x\in V_0(E)}r_{V_0}(g_S)\Phi'_S(1,x)\,r_{V_0}(g^S)\Phi^S(x)
\end{equation*}
is called the inner theta series associated with $A^{(S)}_{\Phi'}$. Note that
here we view $\Phi'_v(1,\cdot)$ as a Schwartz function on $V_1(E_v)$ or $V_0(E_v)$ in two cases. We also set $\theta_{A,0}=0$ if $V_0$ is empty.

These two theta series $\theta_{A,1}(g)$ and $\theta_{A,0}(g)$ can approximate the pseudo-theta series $A^{(S)}_{\Phi'}$ in some sense. In fact, there exists
a finite set $S'$ of non-archimedean places of $F$ and a non-empty open compact
subgroup $K'_S$ of $\UU(\BA_{S'})$ satisfying
\begin{equation*}
 \begin{aligned}
 A_{\Phi'}^{(S)}(g)=&
 \rho_{V,V_1,\infty}(g)\delta(g)^{\dim_E V-\dim_E V_1}\theta_{A,1}(g)\\
 &-\rho_{V,V_0,\infty}(g)\delta(g)^{\dim_E V-\dim_E V_0}\theta_{A,0}(g),
 \end{aligned}
 \qquad g\in K'_{S'}\UU(\BA^{S'}).
\end{equation*}

Let $\ell_v$ be defined by $\chi_v(z)=z^{\ell_v}$ on the norm-one subgroup of $E_v^\times$.

For $g_v=n(b)m(a)[k_1,k_2]$, put
\begin{equation*}
 \rho_v(g_v)=k_1^{(1+\ell_v)/2}k_2^{(1-\ell_v)/2},
 \qquad
 \rho_\infty(g)=\prod_{v\mid\infty}\rho_v(g_v).
\end{equation*}
For a Hermitian space $T$, the compatible splitting character $\chi_T=\chi^{\dim_E T}$ gives
\begin{equation*}
 \mathfrak m_{T,v}=\left(\frac{(\dim_E T)(1+\ell_v)}2,\frac{(\dim_E T)(1-\ell_v)}2\right),
 \qquad
 \rho_{T,\infty}(g)=\rho_\infty(g)^{\dim_E T}.
\end{equation*}
Consequently
\begin{equation*}
 \rho_{V,V_i,\infty}(g)
 =\frac{\rho_{V,\infty}(g)}{\rho_{V_i,\infty}(g)}
 =\rho_\infty(g)^{\dim_E V-\dim_E V_i}.
\end{equation*}
This is the archimedean weight-correction factor in the approximation formula and is compatible with Section \ref{Holomorphic projection}. If $V_0=\emptyset$, the inner-theta term is omitted.

For more details,
see \cite[Section 6.2]{YZ1}.

\subsubsection*{Pseudo-Eisenstein series}
Keep all the notations above, a pseudo-Eisenstein series is a series of the form
\begin{equation*}
    B^{(S)}_\Phi(g)=\sum_{a\in F^\times}B_{a,S}(g)W_a^S(0,g,\Phi^S),\ g\in\UU(\BA).
\end{equation*}
We explain the notations as follows:
\begin{itemize}
    \item $W_a^S(0,g,\Phi^S)=\prod_{w\notin S}W_{a,w}(0,g,\Phi_w)$ is the product of the local Whittaker functions.
    \item $B_{a,S}(g)=\prod_{v\in S}B_{a,v}(g)$ as local product, with
    \begin{equation*}
        B_{\cdot,v}(\cdot):\Fn^\times \times\UU(\Fn)\longrightarrow\CC
    \end{equation*}
    a smooth locally constant function for each $v\in S$. Here the smoothness means the right translation of some open compact subgroup $K_v\subset\UU(\Fn)$ is invariant.
    \item  For any $v\in S$ and $g\in\UU(F_v)$,  the sum is absolutely
    convergent.
\end{itemize}
Note that $B^{(S)}_\Phi(g)$ does not have a constant term.

Recall that in the discussion of Section \ref{Notation}, for a Hermitian space
$(V,q)$ and every $v\in S$, there are one or two Hermitian space over $E_v$
up to isomorphism with the same dimension, namely $(V_v^+,q^+)$
and $(V_v^-,q^-)$ with Hasse invariant 1 or $-1$. Especially when $v$ is split in $E/F$, there is no such $V_v^-$ and we ignore this notation.

The pseudo-Eisenstein series $B^{(S)}(g)$ is called \textit{non-singular} if for every
$v\in S$, there exist $\Phi_v^+\in \mathcal{S}(V_v^+)$ and $\Phi_v^-\in\mathcal{S}(V_v^-)$ such that
\begin{equation*}
    B_{a,v}(1)=W_{a,v}(0,1,\Phi_v^+)+W_{a,v}(0,1,\Phi_v^-),\ \forall a\in\Fn^\times.
\end{equation*}
We take the convention that $W_{a,v}(0,1,\Phi_v^-)=0$ if $V_v^-$ does not exist. Once this is true, we have
\begin{equation*}
    B^{(S)}_\Phi(g)=\sum_{\epsilon:S\rightarrow\{\pm\}}E_*(0,g,\Phi_S^\epsilon\otimes\Phi^S),\ \forall g\in 1_S\UU(\BA^S).
\end{equation*}
Here $\Phi_S^\epsilon=\otimes_{v\in S}\Phi_v^{\epsilon(v)}$. The symbol $V^\epsilon$ denotes the adelic Hermitian space whose $v$-component is $V_v^{\epsilon(v)}$ for $v\in S$ and whose components outside $S$ agree with those of $V$; it is not a tensor product of local vector spaces. The Eisenstein series $E(s,g,\Phi_S^\epsilon\otimes\Phi^S)$ is defined from these adelic Weil data whether $V^\epsilon$ is coherent or incoherent. Coherence is required only when one invokes a theta-integral interpretation through the Siegel--Weil formula.

Meanwhile, if such $\Phi_v^+\in \mathcal{S}(V_v^+)$ and $\Phi_v^-\in\mathcal{S}(V_v^-)$ do not exist, the pseudo-Eisenstein series is called \textit{singular}. Throughout this section, we only treat non-singular pseudo Eisenstein series.

For convenience, denote by
\begin{equation*}
    E_B(g)=E_{B^{(S)}_\Phi}(g)=\sum_{\epsilon:S\rightarrow\{\pm\}}E(0,g,\Phi_S^\epsilon\otimes\Phi^S).
\end{equation*}
This is called the Eisenstein series associated with $B^{(S)}_\Phi$. Although the preimage pair $(\Phi_v^+,\Phi_v^-)$ need not be unique, \cite[Lemma 2.2]{Yuan1} shows that $r(g_v)\Phi_v^+(0)+r(g_v)\Phi_v^-(0)$ is uniquely determined by the nonzero Whittaker coefficients. Hence $E_B(g)$ is independent of the chosen preimages.

Similarly, there is also an approximation property of the non-singular pseudo-Eisenstein series $B^{(S)}_\Phi(g)$. Keep the same notations, we have
\begin{equation*}
 \begin{aligned}
    B^{(S)}_\Phi(g)=E_B(g)
    &+\rho_\infty(g)^{\dim_E V}\delta(g)^{\dim_E V}f_1(g)\\
    &+\rho_\infty(g)^{\dim_E V}\delta(g)^{2-\dim_E V}f_2(g),
 \end{aligned}
 \qquad g\in K'_{S'}\UU(\BA^{S'}).
\end{equation*}
Here $f_1(g)$ and $f_2(g)$ are both automorphic in $g\in\UU(\BA)$, which come from the constant term of the non-singular pseudo-Eisenstein series. Note that both sides are invariant under the action of some non-empty open compact subgroup $K'_S$ of $\UU(\BA_{S'})$, so the approximation actually holds for all $g\in K'_{S'}\UU(\BA^{S'})$.

For this approximation property, see \cite[Section 2.3]{Yuan1}. There is no essential
difference after dropping the similitude $u$ defined there.

\subsubsection*{Key lemma}
The following lemma is a direct corollary of \cite[Lemma 2.3]{Yuan1}. See also \cite[Lemma 6.1(1)]{YZ1}. This lemma will play a fundamental role in the proof of the modular height formula.

\begin{lemma}\label{Key lemma}
Let $\{A_l^{(S_l)}\}_l$ be a finite family of non-singular pseudo-theta series attached to $V_{l,0}\subset V_{l,1}\subset V_l$, and let $\{B_j^{(S'_j)}\}_j$ be a finite family of non-singular pseudo-Eisenstein series. Suppose
\begin{equation*}
 f(g)=\sum_lA_l^{(S_l)}(g)+\sum_jB_j^{(S'_j)}(g)
\end{equation*}
is automorphic. Then
\begin{equation*}
 f(g)=\sum_{l:\,V_{l,1}=V_l}\theta_{A_l,1}(g)+\sum_jE_{B_j}(g).
\end{equation*}
In particular, the degenerate pseudo-theta terms do not appear on the right-hand side.
\end{lemma}
\begin{proof}
This is the Hermitian analogue of \cite[Lemma 2.3]{Yuan1}. The proof uses only the local approximation formulas and automorphy; after fixing compatible splitting characters, the same argument applies to $\UU$. Equivalently, one may restrict to the derived group $\mathrm{SU}(1,1)\simeq\SL_2$ and keep track of the central character.
\end{proof}

\subsubsection*{Local preimages}

Let $v$ be a non-archimedean place of $F$, and let $V_v^+$ and $V_v^-$ be the two local Hermitian spaces of the prescribed dimension and determinant (with the usual convention at split places). Consider
\begin{equation}\label{pseudo Eisenstein preimage}
\mathcal S(V_v^+)\oplus\mathcal S(V_v^-)
\longrightarrow C^\infty(F_v^\times),\qquad
(\Phi^+,\Phi^-)\longmapsto
\bigl(a\mapsto W_{a,v}(0,1,\Phi^+)+W_{a,v}(0,1,\Phi^-)\bigr).
\end{equation}
\begin{lemma}\label{Example of pseudo-Eisenstein series}
Every $\Psi\in\mathcal S(F_v)$ has a preimage under \eqref{pseudo Eisenstein preimage}, which may be chosen so that
\begin{equation}\label{compact open cell origin zero}
 \Phi^+(0)+\Phi^-(0)=0.
\end{equation}
More generally, if a function on $F_v^\times$ lies in the image of \eqref{pseudo Eisenstein preimage} and extends to an element of $\mathcal S(F_v)$, then every preimage satisfies \eqref{compact open cell origin zero}.
\end{lemma}
\begin{proof}
For existence, let $\widehat\Psi$ be the Fourier transform with the convention used in the local Whittaker integral.  In the local degenerate principal series at the Siegel--Weil point, prescribe on the open Bruhat cell
\[
        h(wn(b))=\widehat\Psi(b).
\]
Induction equivariance determines a smooth section there; since $\widehat\Psi$ is compactly supported, it extends across the closed cell with $h(1)=0$.  Fourier inversion gives $W_{a,v}(h)=\Psi(a)$ for all $a\in F_v$.

It remains only to realize $h$ by ordinary Siegel--Weil sections.  In the notation of \cite{GQT}, our group is $G(W_1)$, $d(1)=1$, and for Hermitian dimension $n+1$ the Siegel--Weil point is $s_{n+1,1}=n/2$, which is our $s=0$.  If $n>1$, \cite[Proposition~5.1]{GQT} makes the relevant principal series irreducible, so the nonzero Siegel--Weil image is the whole representation; if $n=1$, \cite[Proposition~5.2(iii)]{GQT} says that the sum of the two opposite-Hasse Siegel--Weil images is the whole principal series.  Thus $h$ is represented by a pair $(\Phi^+,\Phi^-)$ as above, and $h(1)=0$ gives \eqref{compact open cell origin zero}.

For the last assertion, let $(\Phi^+,\Phi^-)$ be any preimage and put $g_b=n(b)m(-b)wn(b)$.  As $v(b)\to-\infty$, one has $g_b\to1$, while
\[
 r(g_b)\Phi^+(0)+r(g_b)\Phi^-(0)
 =\chi_{V_v}(-b)|b|_v^{n+1}
 \bigl(r(wn(b))\Phi^+(0)+r(wn(b))\Phi^-(0)\bigr).
\]
Up to the fixed Fourier-transform normalization, the parenthesized term is the Fourier transform of the Schwartz extension of the coefficient function, hence vanishes for $v(b)$ sufficiently negative.  Passing to the limit gives $\Phi^+(0)+\Phi^-(0)=0$; compare \cite[Remark~2.5]{Yuan1}.
\end{proof}

\begin{remark}
    As we have stated in the proof, there are certain differences between the Hermitian lattice in a Hermitian space and the maximal order in a quaternion algebra, and we will further explore the implications of this difference in our subsequent discussions. Meanwhile, it is worth noticing that the proof of Lemma \ref{Key lemma} does not involve such difference.
\end{remark}

\section{Derivative series}\label{Derivative series}
In this section we study the holomorphic projection of the derivative of a certain mixed Eisenstein-theta series. As introduced in the introduction, we refer to this series as the \textit{derivative series}, which is the analytic side in our final comparison.

\subsection{Derivative series and its decomposition}\label{Derivative series and its decomposition}
In this subsection, we first introduce the derivative series after some preliminary work. Then we apply the local Siegel--Weil formula to obtain a decomposition of such series.

\subsubsection*{Mixed Eisenstein--theta series}
We keep all the notations in Section \ref{Theta}.
For any $\Phi\in\mathcal{S}(\BV)$, we now form a mixed Eisenstein--theta series on $\UU$ associated with the orthogonal decomposition $\BV=\BW\oplus W^\perp(\BA_E)$:
\begin{equation*}
    I(s,g,\Phi)=\sum_{\gamma\in P(F)\backslash\UU(F)}\delta(\gamma g)^s\sum_{x\in W^\perp}r(\gamma g)\Phi(x).
\end{equation*}
In the case of $\Phi=\phi_1\otimes\phi_2$ for $\phi_1\in\mathcal{S}(W^\perp(\BA_E))$ and $\phi_2\in\mathcal{S}(\BW)$, the splitting of the Weil representation implies
\begin{equation*}
    I(s,g,\Phi)=\theta(g,\phi_1)E(s,g,\phi_2).
\end{equation*}

We call $I'(0,g,\Phi)$ the unprojected derivative; the term \emph{derivative series} is reserved for its cuspidal holomorphic projection $\Pr I'(0,g,\Phi)$.

\subsubsection*{Normalization of local Whittaker functions}
Recall our local Whittaker integral
\begin{equation*}
    W_{a,v}(s,g,\phi_{2,v})=\int_{F_v}\delta_v(wn(b)g)^s r(wn(b)g)\phi_{2,v}(0)\psi(-ab)db.
\end{equation*}
Here $\phi_{2,v}\in\mathcal{S}(\BW_v)$. We normalize these local Whittaker functions as
follows. For $a\in\Fn^\times$, denote
\begin{equation*}
    W^\circ_{a,v}(s,g,\phi_{2,v})=\gamma_v^{-1}W_{a,v}(s,g,\phi_{2,v}).
\end{equation*}
Here $\gamma_v$ is the Weil index for $(\BW_v,q_v)$.
For $a=0$, we use the normalization
\begin{equation}\label{Normalized Whittaker function 0}
    W_{0,v}^\circ(s,g,\phi_{2,v})=\gamma_v^{-1}\frac{L(s+1,\eta_v)}{L(s,\eta_v)}|D_v|^{-\frac{1}{2}}|d_v|^{-\frac{1}{2}}W_{0,v}(s,g,\phi_{2,v}).
\end{equation}
Here $D_v$ is the discriminant of $E_v/F_v$, $d_v$ is the different of $F_v/\QQ_p$, and we use $|D_v|=|d_v|=1$ for archimedean place.
Recall that the quadratic character $\eta$ is associated with $E/F$. Note that this normalizing factor
$\frac{L(s+1,\eta_v)}{L(s,\eta_v)}$ has a zero at $s=0$ when $E_v$ is split, and is equal
to $\pi^{-1}$ at $s=0$ if $v$ is archimedean.

Globally, denote
\begin{equation*}
    W^\circ_{a}(s,g,\phi_{2})=\prod_v W^\circ_{a,v}(s,g,\phi_{2,v}).
\end{equation*}
Since the incoherence condition implies $\prod_v \gamma_v=-1$, we have
\begin{equation*}
    W^\circ_{a}(s,g,\phi_{2})=-W_{a}(s,g,\phi_{2}),\quad a\in F^\times.
\end{equation*}
By the functional equation for $L(s,\eta)$, we also have
\begin{equation*}
    W^\circ_{0}(s,g,\phi_{2})=-\frac{L(s+1,\eta)/L(1,\eta)}{L(s,\eta)/L(0,\eta)}W_{0}(s,g,\phi_{2}).
\end{equation*}
Especially for $s=0$, we have
\begin{equation*}
    W^\circ_{0}(0,g,\phi_{2})=-W_{0}(0,g,\phi_{2}).
\end{equation*}

\subsubsection*{Explicit local Siegel--Weil}
Now we present a result without proof. Note that this formula is commonly referred to as ``local Siegel--Weil formula", because local Whittaker functions are the constituents appearing in the expansion of Eisenstein series, while the integration of Schwartz function can be viewed as local theta series. See \cite[Proposition 6.2]{Ich} also.
\begin{proposition}\cite[Proposition 2.7, Proposition 6.1]{YZZ2}\label{local Siegel-Weil}
    \begin{enumerate}
        \item In the sense of analytic continuation for $s\in\CC$,
    \begin{equation*}
        W_{0,v}^\circ(0,g,\phi_{2,v})=r(g)\phi_{2,v}(0).
    \end{equation*}
    Therefore,
    \begin{equation*}
        W_{0}^\circ(0,g,\phi_2)=r(g)\phi_{2}(0).
    \end{equation*}
    Furthermore, the following identity holds at every archimedean place with the standard Gaussian and at all but finitely many non-archimedean places with unramified standard data:
    \begin{equation*}
        W_{0,v}^\circ(s,g,\phi_{2,v})=\delta_v(g)^{-s}r(g)\phi_{2,v}(0).
    \end{equation*}
        \item Assume $a\in F^\times_v$, if $a$ is not represented by $\BW_v$, then $W_{a,v}^\circ(0,g,\phi_{2,v})=0$. If there exists some $x_a\in\BW_v$ satisfying $q(x_a)=a$. Then
    \begin{equation*}
        W_{a,v}^\circ(0,g,\phi_{2,v})=\frac{1}{L(1,\eta_v)}\int_{E_v^1}r(g,h)\phi_{2,v}(x_a)dh.
    \end{equation*}
    Here the integration uses the Haar measure on $E_v^1$ normalized in \ref{Notation}.
    \end{enumerate}
\end{proposition}

A direct corollary is
\begin{equation*}
 E(0,g,\phi_2)=0.
\end{equation*}
Indeed, every nonzero Fourier coefficient vanishes because an element of $F^\times$ fails to be represented at some place of the incoherent Hermitian line $\BW$. The constant term vanishes separately by the normalized functional equation, equivalently by $W_0^\circ(0,g,\phi_2)=-W_0(0,g,\phi_2)$.

Meanwhile, if $a\in F^\times$ can be represented by $\BW(\BA^v)$, i.e, $a$ can be represented
at all places except $v$, then $v$ is nonsplit in $E/F$. Indeed, if
$E_v/F_v$ is split, then $q(E_v)=F_v$ by definition, which implies $a\in F^\times$ can be
represented at every place. This is the reason why the later decomposition is indexed by
the nonsplit place.

We remark that, in the later discussion, we will also use the local Siegel--Weil formula for $\Phi_v\in \mathcal{S}(\BV_v)$. Most of the discussion remains unchanged, except that there are some small modifications at the normalizing factors. See \cite[Proposition 2.9]{YZZ2} for example.

\subsubsection*{A rough decomposition of the derivative}
Clearly, there is a decomposition of our Eisenstein series
\begin{equation*}
    E'(0,g,\phi_2)=E'_0(0,g,\phi_2)-\sum_v\sum_{a\in F^\times} W_{a,v}'(0,g,\phi_{2,v})W_a^v(0,g,\phi_2^v).
\end{equation*}
Following the idea in \cite[Section 6.1.2]{YZZ2}, we have a decomposition of the mixed Eisenstein--theta series
\begin{equation*}
    I'(0,g,\Phi)=-\sum_{v\ \nonsplit}I'(0,g,\Phi)(v)+\theta(g,\phi_1)E'_0(0,g,\phi_2),
\end{equation*}
where for any place $v$ nonsplit in $E$,
\begin{equation*}
    I'(0,g,\Phi)(v)=\theta(g,\phi_1)E'(0,g,\phi_2)(v),
\end{equation*}
\begin{equation*}
    E'(0,g,\phi_2)(v)=\sum_{a\in F^\times}{W_{a,v}^\circ}'(0,g,\phi_{2,v})W_a^{\circ,v}(0,g,\phi_2^v).
\end{equation*}
Note that here we use the equation $W^\circ_{0}(0,g,\phi_{2})=-W_{0}(0,g,\phi_{2})$.
Clearly, the first sum comes from those "non-constant terms", while the second
term comes from the "constant term".

In the next two subsections, we will decompose these two terms separately. The key idea in our decomposition is to use the local Siegel--Weil formula in Proposition \ref{local Siegel-Weil}. Moreover, we will see that many local terms in the decomposition behave like the local Whittaker function.

\subsubsection*{Decomposition of the non-constant term}
Fix a place $v$ of $F$ which is nonsplit in $E$. Denote by $W=W(v)$ the coherent nearby Hermitian space of $\BW$ by changing the local Hasse invariant at $v$. Also denote by $V=V(v)$ the nearby Hermitian space of $\BV$.  From the discussion above, we know immediately that
$q(W)\setminus\{0\}$ is precisely the set of $a\in F^\times$ represented by the corresponding local Hermitian lines at every place, by the local-global principle for the norm form of the quadratic extension $E/F$.

Assume $\Phi_v=\phi_{1,v}\otimes\phi_{2,v}$. Apply the result in Proposition \ref{local Siegel-Weil}, we have
\begin{align*}
    & E'(0,g,\phi_2)(v) \\
    =& \sum_{y_2\in E^1\backslash(W-\{0\})}{W^\circ_{q(y_2),v}}'(0,g,\phi_{2,v})W^{\circ,v}_{q(y_2)}(0,g,\phi_2^v) \\
    =& \frac{1}{L^v(1,\eta)}\sum_{y_2\in E^1\backslash(W-\{0\})}{W^\circ_{q(y_2),v}}'
    (0,g,\phi_{2,v})\int_{E^1(\BA^v)}r(g)\phi_2^v(y_2\tau)d\tau \\
    =& \frac{1}{\vol(E^1_v)L^v(1,\eta)}\sum_{y_2\in E^1\backslash(W-\{0\})}
    \int_{E^1(\BA^v)}{W^\circ_{q(y_2\tau),v}}'(0,g,\phi_{2,v})r(g)\phi_2^v(y_2\tau)d\tau \\
    =& \frac{1}{\vol(E^1_v)L^v(1,\eta)}\int_{E^1\backslash E^1(\BA)}\sum_{y_2\in W-\{0\}}{W^\circ_{q(y_2\tau),v}}'(0,g,\phi_{2,v})r(g)\phi_2^v(y_2\tau)d\tau.
\end{align*}
Therefore, we have the following expression for $I'(0,g,\Phi)(v)$
\begin{equation*}
\begin{aligned}
    I'(0,g,\Phi)(v)=& \frac{1}{\vol(E^1_v)L^v(1,\eta)}\sum_{y_1\in W^\perp}r(g)\phi_1(y_1)\\
    &\cdot\int_{E^1\backslash E^1(\BA)}\sum_{y_2\in W-\{0\}}{W^\circ_{q(y_2\tau),v}}'(0,g,\phi_{2,v})r(g)\phi_2^v(y_2\tau)d\tau
\end{aligned}
\end{equation*}
Combine the two sums for $y_1$ and $y_2$, we have
\begin{equation*}
    \begin{aligned}
    I'(0,g,\Phi)(v)=& \frac{1}{\vol(E^1_v)L^v(1,\eta)}\cdot\int_{E^1\backslash E^1(\BA)}\sum_{y=y_1+y_2\in V,\ y_2\ne 0}\\
    &r(g)\phi_{1,v}(y_1){W^\circ_{q(y_2\tau),v}}'(0,g,\phi_{2,v})r(g)\Phi^v(y\tau)d\tau.
\end{aligned}
\end{equation*}

Now, we define
\begin{equation}\label{k series}
    k_{\Phi_v}(g,y)=\frac{L(1,\eta_v)}{\vol(E_v^1)}r(g)\phi_{1,v}(y_1){W^\circ_{q(y_2),v}}'(0,g,\phi_{2,v}),
\end{equation}
with $y=y_1+y_2\in V_v$, $y_2\ne 0$. We also define
\begin{equation}\label{K series}
   \begin{aligned}
    \mathcal{K}_\Phi^{(v)}(g,\tau)=&\mathcal{K}_{r(\tau)\Phi}^{(v)}(g) \\
    =&\sum_{y\in V-W^\perp} k_{r(\tau)\Phi_v}(g,y)r(g,\tau)\Phi^v(y),
    \end{aligned}
\end{equation}
where $\tau\in E^1(\BA)$ is regarded as an element in $\U(W)(\BA)$ or $\U(V)(\BA)$, and $r(\tau)$ is the Weil representation on $W$ or $V$.

We remind the reader that there is no definition for $\phi_{2,v}(y_2)$, since $\phi_{2,v}$ is a Schwartz function on $\BW_v$, which is not isomorphic to  $W_v$ here by definition. One should be careful about our notations, since we use $q$ to denote the Hermitian form for both $\BW_v$ and $W_v$. Also note that definition \ref{k series} and \ref{K series} can be extended to general Schwartz function $\Phi$ (which are not of the form $\phi_1\otimes\phi_2$) by linearity.

Using these notations, we finally conclude that
\begin{equation}\label{I'(v)}
    \begin{aligned}
        I'(0,g,\Phi)(v)
        =\frac{1}{L(1,\eta)}\int_{E^1\backslash E^1(\BA)}\mathcal{K}_\Phi^{(v)}(g,\tau)d\tau.
    \end{aligned}
\end{equation}
The final expression is very similar to the result in \cite[Proposition 6.5]{YZZ2}.

We also remark that, as in \cite[Lemma 6.6]{YZZ2}, the function $k_{\Phi_v}(g,y)$
behaves like Weil representation under the action of $P(\Fn)$ and $\U(W)(\Fn)$.
See Section \ref{Notation} for the Weil representation in our setting.

\subsubsection*{Decomposition of the constant term}
Following the normalization above and Proposition \ref{local Siegel-Weil}, we obtain
\begin{equation*}
    W'_0(0,g,\phi_2)=-\frac{d}{ds}\Big|_{s=0}(\log\frac{L(s,\eta)}{L(s+1,\eta)})r(g)\phi_2(0)-\sum_v {W^\circ_{0,v}}'(0,g,\phi_{2,v})r(g)\phi_2^v(0).
\end{equation*}
We follow the idea in \cite[Section 7.1]{YZ1} to define a new function
\begin{equation}\label{c sequence}
    c_{\Phi_v}(g,y)=r(g)\phi_{1,v}(y){W^\circ_{0,v}}'(0,g,\phi_{2,v})+\log\delta(g_v)r(g)\Phi_v(y).
\end{equation}
Here $y\in W^\perp$ and $\Phi_v=\phi_{1,v}\otimes\phi_{2,v}$ as above. We emphasize that in this definition, the first Weil representation is on $W^\perp_v$, while the second one is on $\BV_v$. We also define a constant
\begin{equation*}
    c_0=\frac{d}{ds}\Big|_{s=0}(\log\frac{L(s,\eta)}{L(s+1,\eta)}).
\end{equation*}

We remark that $c_{\Phi_v}(g,0)$ is a principal series in the sense that
\begin{equation*}
    c_{\Phi_v}(m(a)n(b)g,0)=\chi_{\BV,v}(a)|a|_{E_v}^{\frac{n+1}{2}}c_{\Phi_v}(g,0),\ a\in E_v^\times,b\in F_v.
\end{equation*}
This follows from our normalization of $W^\circ_{0,v}(s,g,\phi_{2,v})$ and Proposition \ref{local Siegel-Weil}. The exponent $\frac{n+1}{2}$ is half the dimension of $\BV$.

We also introduce a new Eisenstein series from $c_{\Phi_v}(g,y)$
\begin{equation}\label{C sequence}
    C(s,g,\Phi)(v)=\sum_{\gamma\in P(F)\backslash\UU(F)}\delta(\gamma g)^s c_{\Phi_v}(\gamma g,0)r(\gamma g^v)\Phi^v(0),
\end{equation}
\begin{equation*}
    C(s,g,\Phi)=\sum_{v\nmid\infty}C(s,g,\Phi)(v).
\end{equation*}
In fact, we know that the second summation has only finitely many nonzero terms since $c_{\Phi_v}(g,0)=0$ for all but finitely many $v$. See Proposition \ref{local Siegel-Weil}. This Eisenstein sequence will be used to understand the absolute constant term of $I'(0,g,\Phi)$.

Let us return to our main discussion. Combine all the discussion, we have the decomposition
\begin{equation} \label{Decomposition of I'}
    \begin{aligned}
    I'(0,g,\Phi)=&-\sum_{v\ \nonsplit}I'(0,g,\Phi)(v)-c_0\sum_{y\in W^\perp}r(g)\Phi(y)\\
    &-\sum_v\sum_{y\in W^\perp} c_{\Phi_v}(g,y)r(g)\Phi^v(y)+2\log\delta(g)\sum_{y\in W^\perp}r(g)\Phi(y).
    \end{aligned}
\end{equation}
One remark is that from the functional equation of the complete L-function $L(s,\eta)$, we have
\begin{equation}\label{c0}
    c_0=2\frac{L'(0,\eta)}{L(0,\eta)}+\log\lvert d_E/d_F\rvert.
\end{equation}

\subsection{Holomorphic projection}\label{Holomorphic projection}
In this subsection we give an explicit expression of the holomorphic projection of $I'(0,g,\Phi)$, which we denote by $\Pr I'(0,g,\Phi)$. Just like in the references, we need to introduce a concept called \textit{quasi-holomorphic projection}, and then utilize this concept to carry out specific calculations.

\subsubsection*{Basic definition}

For each archimedean place $v$, write $\chi_v(z)=z^{\ell_v}$ on the norm-one subgroup of $E_v^\times$. Since $\chi_{\BV,v}=\chi_v^{n+1}$, the weight of $I'(0,g,\Phi)$ is $\mathfrak m=(\mathfrak m_v)_{v\mid\infty}$ with
\begin{equation*}
 \mathfrak m_v=
 \left(\frac{(n+1)(1+\ell_v)}2,
       \frac{(n+1)(1-\ell_v)}2\right).
\end{equation*}
For the fixed finite $K$-type of $I'(0,g,\Phi)$, its restriction to the adelic center has a finite character decomposition
\begin{equation*}
 I'(0,g,\Phi)=\sum_\omega I'(0,g,\Phi)_\omega.
\end{equation*}
Let $\mathcal A_0^{(\mathfrak m)}(\UU(\BA),\omega)$ be the space of holomorphic cusp forms of central character $\omega$ and weight $\mathfrak m$. For an automorphic form $f$ of moderate growth in the $\omega$-component, its cuspidal holomorphic projection $\Pr(f)$ is characterized by
\begin{equation*}
 \langle\Pr(f),h\rangle_{\mathrm{Pet}}=\langle f,h\rangle_{\mathrm{Pet}}
 \qquad
 \text{for all }h\in\mathcal A_0^{(\mathfrak m)}(\UU(\BA),\omega).
\end{equation*}
The pairing converges against cusp forms by their rapid decay, and we apply this definition componentwise. Thus $\Pr$ is not being asserted to be a Hilbert-space orthogonal projection on the unrestricted space of all automorphic forms.

For $t\in F^\times$, put $\psi_t(b)=\psi(tb)$, and let $f_t$ denote the corresponding Whittaker coefficient. Define
\begin{equation*}
 f_{t,s}(g)=(4\pi)^{[F:\QQ]n}\Gamma(n)^{-[F:\QQ]}
 W_t^{(\mathfrak m)}(g_\infty)
 \int_{Z(F_\infty)N(F_\infty)\backslash\UU(F_\infty)}
 \delta(h)^s f_t(gh)\overline{W_t^{(\mathfrak m)}(h)}\,dh.
\end{equation*}
The holomorphic Whittaker function is normalized by
\begin{equation*}
 \int_{Z(F_v)N(F_v)\backslash\UU(F_v)}
 |W_{t,v}^{(\mathfrak m)}(h)|^2dh=(4\pi)^{-n}\Gamma(n)
 \qquad(v\mid\infty).
\end{equation*}
If $f_{t,s}$ is meromorphic at $s=0$, let
\begin{equation*}
 \Pr'(f)_t(g)=\widetilde{\lim}_{s\to0}f_{t,s}(g),
 \qquad
 \Pr'(f)(g)=\sum_{t\in F_{>0}}\Pr'(f)_t(g),
\end{equation*}
where the quasi-limit is the constant Laurent coefficient. Only totally positive $t$ occur in a holomorphic Fourier expansion.

\begin{proposition}\label{Growth condition}
Assume that there is an $\epsilon>0$ such that, for every archimedean place $v$ and with the remaining components fixed,
\begin{equation*}
 f(m(a_v)g)=O_g\bigl(|a_v|_{E_v}^{(n+1)/2-\epsilon}\bigr)
 \qquad\bigl(|a_v|_{E_v}\to\infty\bigr).
\end{equation*}
Then $\Pr(f)=\Pr'(f)$.
\end{proposition}
This is the growth criterion of \cite[Proposition 6.12]{YZZ2} in the present normalization; see also the unitary formulation in \cite[Section 3.3]{Qiu}.

However, our series $I'(0,g,\Phi)$ does not satisfy such growth condition. Thus, we need to track the growth carefully, since these extra terms contributed by the growth will imply our final result.

\subsubsection*{Absolute constant term}
In order to apply Proposition \ref{Growth condition}, the key point is to understand the absolute constant term
\begin{equation}\label{I00 sequence}
    I_{00}(s,g,\Phi)=\theta_0(g,\phi_1)E_0(s,g,\phi_2).
\end{equation}
Recall that we assume $\Phi=\phi_1\otimes\phi_2$. Thus by definition, it follows that
\begin{equation*}
    I_{00}(s,g,\Phi)=\delta(g)^s r(g)\Phi(0)+r(g)\phi_1(0)W_0(s,g,\phi_2),
\end{equation*}
which is a sum of two principal series. In fact, for any $g'=\matrixx{a}{b}{}{\bar{a}^{-1}}$, we have
\begin{equation*}
    \delta(g'g)^s r(g'g)\Phi(0)=\chi_\BV(a)\lvert a\rvert^{\frac{s}{2}+\frac{n+1}{2}}_{\BA_E}\delta(g)^sr(g)\Phi(0)
\end{equation*}
and
\begin{equation*}
    r(g'g)\phi_1(0)W_0(s,g'g,\phi_2)=\chi_\BV(a)\lvert a\rvert^{-\frac{s}{2}+\frac{n+1}{2}}_{\BA_E} r(g)\phi_1(0)W_0(s,g,\phi_2).
\end{equation*}
Here $|c|_{\BA_E}=\prod_v|c_v|_{E_v}$ for $c\in\BA_E^\times$.

Now we define $\mathcal{J}(s,g,\Phi)$ to be the Eisenstein series formed by $I_{00}(s,g,\Phi)$, i.e.,
\begin{equation}\label{J sequence}
    \mathcal{J}(s,g,\Phi)=\sum_{\gamma\in P(F)\backslash\UU(F)} I_{00}(s,\gamma g,\Phi).
\end{equation}

The proof of \cite[Lemma 6.13]{YZZ2} is formal after fixing the central character and replacing its diagonal element by $m(a)$; it gives that the difference
\begin{equation*}
    I'(0,g,\Phi)-I'_{00}(0,g,\Phi)
\end{equation*}
satisfies the growth condition. Note that in the original statement there is a character $\chi$ and the series all have weight $2$, the same estimates apply to the present weights.

Similarly, we can prove that
\begin{equation*}
    \mathcal{J}'(0,g,\Phi)-I'_{00}(0,g,\Phi)
\end{equation*}
satisfies the growth condition. Thus, we conclude that
\begin{equation*}
    I'(0,g,\Phi)-\mathcal{J}'(0,g,\Phi)
\end{equation*}
satisfies the growth condition. Combine with Proposition \ref{Growth condition}, we have
\begin{equation*}
    \Pr I'(0,g,\Phi)=\Pr' I'(0,g,\Phi)-\Pr'\mathcal{J}'(0,g,\Phi).
\end{equation*}

\subsubsection*{Projection of $I'$-part}
The exclusion of $0$ from the sums below is essential in the logarithmic term.
Recall our decomposition of $I'(0,g,\Phi)$ in \ref{Decomposition of I'}. We have the following result concerning the term $\Pr' I'(0,g,\Phi)$.
\begin{theorem}\label{Projection of I'}
    Assume that $\Phi$ is standard at infinity as in \ref{Notation}. Then we have
    \begin{equation*}
        \begin{aligned}
            \Pr' I'(0,g,\Phi)= &-\sum_{v|\infty}\overline{I'}(0,g,\Phi)(v)-\sum_{v\nmid\infty\ \nonsplit}I'(0,g,\Phi)(v)\\
            &-c_1\sum_{y\in W^\perp\setminus\{0\}}r(g)\Phi(y)-\sum_{v\nmid\infty}\sum_{y\in W^\perp\setminus\{0\}}c_{\Phi_v}(g,y)r(g)\Phi^v(y) \\
            &+\sum_{y\in W^\perp\setminus\{0\}}(2\log\delta_f(g_f)+\log\lvert q(y)\rvert_f)r(g)\Phi(y).
        \end{aligned}
    \end{equation*}
    The notations are explained as follows.

    \begin{enumerate}[(1)]
        \item For any archimedean place $v$,
        \begin{equation*}
            \overline{I'}(0,g,\Phi)(v)=\frac{1}{L(1,\eta)}\int_{E^1\backslash E^1(\BA)}\overline{\mathcal{K}}_\Phi^{(v)}(g,\tau)d\tau,
        \end{equation*}
        \begin{equation*}
            \overline{\mathcal{K}}_\Phi^{(v)}(g,\tau)=\wt{\lim_{s\rightarrow0}}\sum_{y\in V-W^\perp}k_{v,s}(y)r(g,\tau)\Phi(y),
        \end{equation*}
        \begin{equation*}
            k_{v,s}(y)=\frac{\Gamma(s+n)}{2(4\pi)^s\Gamma(n)}\int_1^\infty \frac{1}{t(1-q(y_2)t)^{s+n}}dt.
        \end{equation*}
        \item For $I'(0,g,\Phi)(v)$ with $v$ a non-archimedean place, the definition is the same as the one in Section \ref{Derivative series and its decomposition}.
        \item The constant
        \begin{equation}\label{c1}
            c_1=2\frac{L_f'(0,\eta)}{L_f(0,\eta)}+\log\lvert d_E/d_F\rvert.
        \end{equation}
        Here $L_f$ means the finite part of $L$-function.
        \item The definition of $c_{\Phi_v}(g,y)$ is in \ref{c sequence}.
        \item $\lvert\cdot\rvert_f$ means the norm with respect to $\Af$.
    \end{enumerate}
\end{theorem}
\begin{proof}
    The proof of this theorem follows from the one in \cite[Theorem 7.2]{YZ1}. Especially, to assure the readers that the proof of this theorem aligns with the situation of the quaternionic Shimura curve, we will check the explicit expression of $k_{v,s}(y)$ in Lemma \ref{Explicit archimedean k series}.

    It is worth noting a slight difference in the expressions, where we use $q(y_2)$ in our formulation, while in the corresponding part of \cite{YZZ2,YZ1,Yuan1}, they use $\lambda(y)=q(y_2)/q(y)$ instead. The main reason for this phenomenon lies in the different group actions in the two cases. Simply put, the action given by the quaternion algebra itself does not preserve the reduced norm, while the action of the unitary group on the Hermitian space does preserve the Hermitian norm. Thus, there is no fundamental difference between the expressions in the two cases.
\end{proof}

\subsubsection*{Projection of $J'$-part}
Now we compute the term $\Pr'\mathcal{J}'(0,g,\Phi)$. Recall the definition \ref{c sequence} and \ref{I00 sequence}, it follows that
\begin{equation*}
    \begin{aligned}
        &I_{00}'(0,g,\Phi) \\
        =&\log\delta(g)\cdot r(g)\Phi(0)-c_0r(g)\Phi(0)-\sum_v r(g^v)\phi^v(0)\cdot r(g_v)\phi_{1,v}(0){W_{0,v}^\circ}'(0,g_v,\phi_{2,v})\\
        =&2\log\delta(g)\cdot r(g)\Phi(0)-c_0r(g)\Phi(0)-\sum_v r(g^v)\phi^v(0)c_{\Phi_v}(g,0).
    \end{aligned}
\end{equation*}
Combined with definition \ref{J sequence}, we know immediately that
\begin{equation*}
    \mathcal{J}'(0,g,\Phi)=2E'(0,g,\Phi)-c_0 E(0,g,\Phi)-C(0,g,\Phi).
\end{equation*}
We remind the reader that here the Eisenstein series is associated with $\BV$. Apply the formula for $\Pr'$, we have the following expression.

\begin{theorem}\label{Projection of J'}
    Assume that $\Phi$ is standard at infinity. Then
    \begin{equation*}
        \Pr'\mathcal{J}'(0,g,\Phi)=(-c_0+2c_3[F:\QQ])E_*(0,g,\Phi)-C_*(0,g,\Phi)+2\sum_{v\nmid\infty}E'(0,g,\Phi)(v).
    \end{equation*}
    Here $E_*$ and $C_*$ are the non-constant parts of the Eisenstein series $E$ and $C$, $c_3$ is a constant defined as
    \begin{equation}\label{c3}
        c_3=\sum_{i=1}^{n-1}\frac{1}{2i}-\frac{1}{2n}-\frac{\log4}{2},
    \end{equation}
    which is explicitly computed in Lemma \ref{Explicit archimedean}.
    For any non-archimedean $v$, similar to the decomposition of $E'(0,g,\phi_2)$ in
    Section \ref{Derivative series and its decomposition}, we have
    \begin{equation*}
        E'(0,g,\Phi)(v)=\sum_{a\in F^\times}W_a^v(0,g,\Phi^v)\Big(W_{a,v}'(0,g,\Phi_v)-\frac{1}{2}\log|a|_v\cdot W_{a,v}(0,g,\Phi_v)\Big).
    \end{equation*}
\end{theorem}
\begin{proof}
    The proof is parallel to the proof of Theorem \ref{Projection of I'} and \cite[Proposition 3.2]{Yuan1}. The only difference is to compute the constant $c_3$. See \ref{Explicit archimedean}. Thus, we again omit the proof. See also
    \cite[(3.31)]{Qiu}.
\end{proof}

Moreover, as we remarked at the end of Section \ref{Theta}, for $t\in F_{>0}$ denote by $\mathcal J_t(s,g,\Phi)$ the $\psi_t$-Whittaker function of $\mathcal J(s,g,\Phi)$, with analogous notation for $E_t$ and $C_t$. Then Theorem \ref{Projection of J'} remains valid after replacing every term by its $\psi_t$-component.

Combining Theorem \ref{Projection of I'} and Theorem \ref{Projection of J'}, we obtain a decomposition of the derivative series $\Pr I'(0,g,\Phi)$ as a linear combination of pseudo-theta series and pseudo-Eisenstein series. This proves the first part of our main Theorem \ref{main theorem of derivative series}.

\section{Explicit local derivatives}\label{Explicit derivative}
In this section we do some explicit computation under the assumptions we made. The results of these calculations will provide all the information for the analytic side. Specifically, we compute the Eisenstein series and Whittaker function for both $\Phi\in\mathcal{S}(\BV)$ and $\phi_2\in\mathcal{S}(\BW)$ respectively, which also means that we obtain an explicit expression for $\Pr' I'$ and $\Pr' J'$ in Theorems \ref{Projection of I'} and \ref{Projection of J'}. The computation of terms in $\Pr'I'$ is closely related to the computation in \cite[Section 6]{YZZ2} and \cite[Section 7.3]{YZ1}, while the computation of terms in $\Pr'J'$ is a high-dimensional generalization of \cite[Section 3.3]{Yuan1}, which means we need to introduce some additional computational techniques.

All explicit local formulas below are evaluated at $g=1$, which is the value required in the final comparison. The covariance formulas in \eqref{action on Whittaker function}, together with Iwasawa decomposition and the action of the compact factor on the chosen Schwartz function, recover the values on the relevant Iwasawa cells. Varying the Fourier index $a$ alone does not determine the value for an arbitrary $g$ unless its compact factor fixes $\Phi_v$.

\subsection{Choice of the Schwartz function}\label{Choice of Schwartz function}
In this subsection, we introduce some notations and conditions. Some of the assumptions
are already mentioned before, while the other assumptions are restrictive, and will simplify our further computations.

\subsubsection*{Basic notations}
For convenience, we introduce the following notations:
\begin{itemize}
    \item Denote by $\Sigma_\ram$ the set of places in $F$ that are ramified in $E/F$;
    \item Denote by $\Sigma_\Has$ the sets of places in $F$, such that $\epsilon(\BV_v,q)=-1$.
\end{itemize}
In our explicit computation, we always assume that $\Sigma_\Has\subset\Sigma_\ram$, and $|d_v|=1$ if $|D_v|\ne 1$. Here $d_v\in F_v$ is the local different of $F_v$ over $\QQ_p$ and $D_v\in F_v$ is the discriminant of $E_v/F_v$.

Recall that in definition \ref{Choice of lattice}, we fixed a lattice
\begin{equation*}
    \Lambda=\prod_v \Lambda_v\subset\BV_f,
\end{equation*}
such that for $2|n$, $\Lambda_v$ is self-dual at unramified places and almost $\varpi_{E_v}$-modular at ramified places, while for $2\nmid n$, $\Lambda_v$ is self-dual at unramified places and $\varpi_{E_v}$-modular at ramified places. For each finite place $v$, we also fix a vector $e_v\in\BW_v$ such that $q(e_v)=d_{\BV,v}\in\mathcal{O}_{F_v}^\times$. Especially when $v$ is unramified in $E/F$, there exists a decomposition
\begin{equation*}
    \Lambda_v=\Lambda^\perp_{v}\oplus \Lambda_{\BW_v}
\end{equation*}
with $\Lambda_{\BW_v}$ generated by $e_v$.

Moreover, let $U_\Lambda\subset\U(\BV_f)$ be the stabilizer of $\Lambda$, which is an open compact subgroup of $\U(\BV_f)$. Here the subscript $f$ means the product of all finite places.

\subsubsection*{Restrictive conditions}
We now fix the Schwartz datum used in all explicit calculations. Let $\Phi=\otimes_v\Phi_v\in\mathcal S(\BV)$, and take
\begin{enumerate}
 \item the standard Gaussian at every archimedean place;
 \item $\Phi_v=1_{\Lambda_v}$ at every finite place, including the ramified places.
\end{enumerate}
At a bad place, ``standard'' therefore means lattice-standard; it does not mean self-dual under Fourier transform.
Note that if $v$ is unramified in $E/F$, $\Phi_v$ can be written as $\phi_{1,v}\otimes\phi_{2,v}$, where $\phi_1=1_{\Lambda_v^\perp}\in\mathcal{S}(W^\perp(\BA))$ and $\phi_2=1_{\Lambda_{\BW_v}}\in\mathcal{S}(\BW)$. If $v$ is ramified in $E/F$, there is no such decomposition.

Under the natural action of $\U(\BV_f)$ on $\mathcal S(\BV_f)$, the finite part of $\Phi$ is invariant under the lattice stabilizer $U_\Lambda$. There are not separate left and right actions on this Schwartz space. Unlike \cite[Section 3.2]{Yuan1}, we impose no degeneracy condition at the bad places.

\subsection{Explicit computations of I-part: archimedean places}
The following lemma gives the explicit result of $\overline{I'}(0,g,\Phi)(v)$ for archimedean place $v$, which is mentioned in Theorem \ref{Projection of I'}. This lemma is a modification of \cite[Proposition 6.11 and Proposition 6.15]{YZZ2}. The result involves the real exponential integral on the negative axis,
\begin{equation*}
 \Ei(x)=\int_{-\infty}^x\frac{e^t}{t}\,dt
 =\gamma+\log(-x)+\int_0^x\frac{e^t-1}{t}\,dt,
 \qquad x<0.
\end{equation*}
This is the only range used below, since $q(y_2)<0$.

\begin{lemma}\label{Explicit archimedean k series}
    Let $v$ be an archimedean place, and in this case, $W_v$ is a negative definite Hermitian space of dimension 1, where $W$ is the nearby coherent
    Hermitian space of $\BW$ with respect to $v$.
    \begin{enumerate}
        \item We have
        \begin{equation*}
           k_{\Phi_v}(g,y)=-\frac{1}{2}|a|^{n+1}e^{-2\pi q(y) |a|^2}e^{2\pi iq(y)b}k_1^{\mm_{v,1}}k_2^{\mm_{v,2}}\Ei(4\pi q(y_2)|a|^2).
        \end{equation*}
        for any
        \begin{equation*}
            g=\matrixx{a}{}{}{\bar{a}^{-1}}\matrixx{1}{b}{}{1}[k_1,k_2]\in\UU(\Fn)
        \end{equation*}
        in the form of the Iwasawa decomposition. Here $|\cdot|$ is the usual norm on $\CC$, and $y=y_1+y_2$ is an element
        in $V_v=W_v^\perp\oplus W_v$ with $y_2\ne 0$.
        \item The term $k_{v,s}(y)$ used in $\overline{I'}(0,g,\Phi)(v)$ is defined by
        \begin{equation*}
            k_{v,s}(y)=\frac{\Gamma(s+n)}{2(4\pi)^s\Gamma(n)}\int_1^\infty \frac{1}{t(1-q(y_2)t)^{s+n}}dt.
        \end{equation*}
    \end{enumerate}
\end{lemma}

\begin{proof}
    We refer to \cite[Proposition 2.11]{YZZ2} for general result. To prove our
    first result, it is sufficient to check the formula in the case $g=1$:
    \begin{equation*}
        k_{\Phi_v}(1,y)=-\frac{1}{2}e^{-2\pi q(y)}\Ei(4\pi q(y_2)).
    \end{equation*}
    Then the general case is obtained by the behavior of $k_{\Phi_v}(g,y)$
    under the action of $P(\Fn)$ and $\U(W)(\Fn)$, which is the same as the behavior of Weil representation mentioned before. The character $\chi$ is trivial here. Also note that we have
    \begin{equation*}
        r([k_1,k_2])\Phi_v=k_1^{\mm_{v,1}}k_2^{\mm_{v,2}}\Phi_v.
    \end{equation*}
    Here $[k_1,k_2]$ is defined in Section \ref{Notation}, and $\mm_v=(\mm_{v,1},\mm_{v,2})$ is defined in Section \ref{Holomorphic projection}.

    Following the definition of $k_{\Phi_v}(1,y)$, recall that $\vol(E^1_v)=2$ and $L(1,\eta_v)=\pi^{-1}$, it amounts to show that for
    any $a<0$,
    \begin{equation*}
        {W_{a,v}^\circ}'(0,1,\phi_{2,v})=-\pi e^{-2\pi a}\Ei(4\pi a).
    \end{equation*}
    This is just the case of $d=2$ in \cite[Proposition 2.11 (3)]{YZZ2}.

    Note that here we only need to consider the case $a<0$ since $W=W(v)$ is the
    coherent nearby Hermitian space of $\BW$ by changing the local Hasse invariant at $v$, which means $(W_v,q)$ is negative definite. This also implies the signature here is $-1$.

    To compute the quasi-holomorphic projection, it is sufficient to consider
    \begin{equation*}
        \wt{k}_{\Phi_v,s}(g,y)=(4\pi)^{n}\Gamma(n)^{-1}W^{(\mm_v)}(g)\int_{Z(\RR)N(\RR)\backslash\UU(\RR)}\delta(h)^s k_{\Phi_v}(gh,y)
    \overline{W^{(\mm_v)}(h)}dh.
    \end{equation*}
    By the equivariance under the diagonal torus, the calculation of the positive Fourier coefficient may be normalized to $q(y)=1$; the general coefficient is recovered by scaling. A direct computation then implies
    \begin{equation*}
        \begin{aligned}
            \wt{k}_{\Phi_v,s}(g,y)=&(4\pi)^{n}\Gamma(n)^{-1}W^{(\mm_v)}(g)\int_{t>0}t^{s+\frac{n-1}{2}} e^{-2\pi t}k_{\Phi_v}(m(\sqrt{t}),y)\frac{dt}{t} \\
            =&-\frac{1}{2}(4\pi)^{n}\Gamma(n)^{-1}W^{(\mm_v)}(g)\int_{t>0}t^{s+n} e^{-4\pi t}\Ei(4\pi q(y_2)t)\frac{dt}{t} \\
            =&\frac{1}{2}(4\pi)^{n}\Gamma(n)^{-1}W^{(\mm_v)}(g)\int_{t>0}t^{s+n} e^{-4\pi t}\int_1^\infty e^{4\pi q(y_2)t t'}\frac{dt'}{t'}\frac{dt}{t} \\
            =&\frac{1}{2}(4\pi)^{n}\Gamma(n)^{-1}W^{(\mm_v)}(g)\int_1^\infty t'^{-1}\int_{0}^\infty t^{s+n} e^{-2\pi t(-2q(y_2)t'+2)}\frac{dt}{t}dt' \\
            =&\frac{1}{2}(4\pi)^{-s}\frac{\Gamma(s+n)}{\Gamma(n)}W^{(\mm_v)}(g)\int_1^\infty  \frac{1}{t(-q(y_2)t+1)^{s+n}}dt\\
        \end{aligned}
    \end{equation*}
    Hence, we conclude our second result.
\end{proof}

    We remark that, by the definition of $W_v$ here, $q(y_2)<0$ for any nonzero $y_2\in W_v$.

\subsection{Explicit computations of I-part: non-archimedean places}
Now, assume $v$ is a non-archimedean place of $F$. Since we only have the term $I'(0,g,\Phi)(v)$ for nonsplit places, we further assume $v$ is either inert or ramified in $E/F$. The following lemma \ref{Explicit nonarchimedean k series} is similar to \cite[Corollary 6.8]{YZZ2} and \cite[Lemma 7.4]{YZ1}, which gives an explicit expression of $k_{\Phi_v}(1,y)$ defined in \ref{k series}.

We first emphasize some notations here. Recall that under our notation convention, we denote by $V=V(v)$ the nearby coherent Hermitian space of $\BV$ with respect to a non-archimedean place $v$, and denote by $W=W(v)$ the Hermitian space similarly. Note that here
\begin{equation*}
    \BW^\perp_v\cong W^\perp_v,\  \BW_v\ncong W_v.
\end{equation*}

Recall that $\Lambda_{\BW_v}=\mathcal{O}_{E_v}e_v\subset\BW_v$. Following \ref{decomposition of inert Hermitian lattice} and \ref{decomposition of ramified Hermitian lattice}, we again remind the reader that when $v$ is ramified in $E/F$, $\Lambda_{\BW_v}$ is \textit{not} a direct summand of $\Lambda_v$. Denote by $\Lambda_{W_v}\subset W_v$ a nearby lattice of $\Lambda_{\BW_v}$ in the nearby space, while the corresponding lattice of $\Lambda^\perp_v$ is isomorphic to itself, so we still use the notation $\Lambda^\perp_v$.

Let us briefly explain the definition of this nearby lattice $\Lambda_{W_v}$. If $v$ is inert in $E/F$, since we assume $\BW_v$ is a split Hermitian space, i.e., $\epsilon(\BW_v,q)=1$, then $W_v$ is nonsplit. In this case, fix a vector $\mathfrak{e}_v\in W_v$ such that $v(q(\mathfrak{e}_v))=1$, and we define $\Lambda_{W_v}=\mathcal{O}_{E_v}\mathfrak{e}_v$. If $v$ is ramified in $E/F$, regardless of whether $\BW_v$ is split or not, we can always find a vector in $\mathfrak{e}\in W_v$ such that $v(q(\mathfrak{e}_v))=0$. Then we also define $\Lambda_{W_v}=\mathcal{O}_{E_v}\mathfrak{e}_v$.

\begin{lemma}\label{Explicit nonarchimedean k series}
    Let $y=y_1+y_2$ be an element in $V_v=W^\perp_v\oplus W_v$ with $y_2\ne 0$.
    \begin{enumerate}
        \item Let $v$ be a non-archimedean place inert in $E/F$. Then the difference
        \begin{equation*}
            k_{\Phi_v}(1,y)-\phi_{1,v}(y_1)\cdot 1_{\Lambda_{W_v}}(y_2)\cdot\frac{1}{2}\big(v(q(y_2))+1\big)\log N_v
        \end{equation*}
        extends to a Schwartz function on $V_v$ whose restriction to $\Lambda_v^\perp$ is equal to
        \begin{equation*}
            \frac{|d_v|-1}{(1-N_v)(1+N_v^{-1})}\cdot\phi_{1,v}(y_1)\cdot\log N_v.
        \end{equation*}
        \item Let $v$ be a non-archimedean place which is ramified in $E/F$. Then the difference
        \begin{equation*}
             k_{\Phi_v}(1,y)-\Phi_{v}(y_1)\cdot 1_{\Lambda_{W_v}}(y_2)\cdot\frac{1}{2}\big(v(q(y_2))+1\big)\log N_v
        \end{equation*}
        extends to a Schwartz function on $V_v$ whose restriction to $\Lambda_v^\perp$ is $0$.
    \end{enumerate}
\end{lemma}

\begin{proof}
    By the compatible local choice above and the lattice classification in \cite[Lemma 5.1.1]{Qiu},
\begin{equation*}
 \Lambda_v=\Lambda_{1,v}\oplus\Lambda_{1,v}^\perp
\end{equation*}
with $W_v\subset V_{1,v}$, and the lattice-standard Schwartz function factors accordingly. The nontrivial local Whittaker calculation therefore occurs on the rank-two factor $V_{1,v}$; the complementary factor contributes its value at the origin. Thus the calculation reduces to the local case $n=1$. This is a reduction using the chosen orthogonal factorization, not a statement that an arbitrary lattice stabilizer acts transitively on all vectors.

    Thanks to this assumption, we can now directly utilize the results from \cite[Lemma 7.4]{YZ1}. The proof of the inert case is entirely the same, since we have a decomposition
    \begin{equation*}
        \Lambda_{v}=\Lambda_{\BW_v}\oplus \mathcal{O}_{E_v}e_v'\subset\BV_{v},
    \end{equation*}
    where $e_v'\in\BV_{v}$ satisfies $q(e_v')\in\mathcal{O}_{F_v}^\times$. Then our result of inert case becomes a special case of \cite[Lemma 7.4(1)]{YZ1} when taking $|q(\mathfrak{j}_v)|=1$. If readers still have doubts about how to directly use the proof from \cite[Lemma 7.4]{YZ1}, we refer to the proof of Lemma \ref{Explicit nonarchimedean}. There, we will provide a detailed proof.

    It remains to prove the ramified case. Put
    \begin{equation*}
        A_v=\mathcal O_{E_v}e_v\oplus\mathcal O_{E_v}e'_v.
    \end{equation*}
    Then \eqref{decomposition of ramified Hermitian lattice} reads
    \begin{equation*}
        \varpi_{E_v}A_v\subset\Lambda_{1,v}\subset A_v.
    \end{equation*}

Let
\begin{equation*}
 \Lambda_v'=A_v\oplus\Lambda_{1,v}^\perp,
 \qquad
 \Phi_v'=1_{\Lambda_v'}.
\end{equation*}
The rank-two lattice $A_v$ is decomposable with respect to
$V_{1,v}=W_v^\perp\oplus W_v$, so the ramified formula of
\cite[Lemma 7.4(2)]{YZ1} applies directly to $\Phi_v'$. To pass from
$\Phi_v'$ to the actual lattice-standard function $\Phi_v=1_{\Lambda_v}$,
one must use the two colength-one inclusions
\begin{equation*}
 \varpi_{E_v}A_v\subset\Lambda_{1,v}\subset A_v;
\end{equation*}
bare linearity does not by itself identify the required extension at the origin.
The comparison is exactly the argument of \cite[Lemma 5.2.3]{Qiu}:
one restricts the extension for $1_{A_v}$ to the smaller support, separates
the two boundary pieces of $A_v/\Lambda_{1,v}$, and applies
\cite[Lemma 3.1.2]{Qiu} to the resulting Schwartz functions that vanish at
the origin. This proves that the displayed difference in the statement
extends to a Schwartz function on $V_v$. Under our standing assumptions
$|d_v|=1$ and $v(D_v)=1$, the correction term in
\cite[Lemma 7.4(2)]{YZ1} has zero restriction to $\Lambda_v^\perp$, which
gives the asserted value there. The argument is valid whether the full
ambient lattice is $\varpi_{E_v}$-modular or almost
$\varpi_{E_v}$-modular, because the comparison takes place on the chosen
rank-two modular summand $\Lambda_{1,v}$.

\end{proof}

\begin{remark}\label{quaternion and Hermitian remark}
    In fact, there is a miracle correspondence between quaternion algebra and Hermitian space of dimension 2, which was first suggested by Deligne. Using this correspondence, the above proof becomes very natural. Moreover, this identification has many applications in work on special value formulas, and more broadly in problems related to Shimura varieties in general. For instance, \cite{How} use this correspondence to derive the arithmetic volume of unitary Shimura curve from the one of quaternionic Shimura curve. In our series of papers, we also use this correspondence in the second paper \cite{Guo1}.
\end{remark}

The remaining three terms in $\Pr'I'$,
\begin{equation*}
\begin{aligned}
 &-c_1\sum_{y\in W^\perp\setminus\{0\}}r(g)\Phi(y)
 -\sum_{v\nmid\infty}\sum_{y\in W^\perp\setminus\{0\}}
   c_{\Phi_v}(g,y)r(g)\Phi^v(y)\\
 &\hspace{3em}
 +\sum_{y\in W^\perp\setminus\{0\}}
   \bigl(2\log\delta_f(g_f)+\log|q(y)|_f\bigr)r(g)\Phi(y),
\end{aligned}
\end{equation*}
are degenerate pseudo-theta series supported on the proper subspace $W^\perp$. They are not individually zero and need not have zero Fourier constant term. Rather, in the automorphic combination to which Lemma \ref{Key lemma} is applied, degenerate pseudo-theta terms do not occur on the theta-series side and hence do not enter the final constant-term identity.

\subsection{Explicit computations of J-part: archimedean places}
The following lemma computes the constant $c_3$ used in Theorem \ref{Projection of J'}. This lemma is a higher-dimensional generalization of \cite[Lemma 3.3]{Yuan1}.
\begin{lemma}\label{Explicit archimedean}
    Let $v$ be an archimedean place.
    \begin{enumerate}
        \item For any $a\in\RR^\times$ (the constant term $a=0$ is treated separately),
        \begin{equation*}
    W_{a,v}(0,1,\Phi_v)=\left\{
    \begin{aligned}
        \nonumber
        &e^{\frac{2\pi i(1+n)}{4}}\frac{(2\pi)^{1+n}}{\Gamma(1+n)}a^{n}e^{-2\pi a}\quad a>0,\\
        &0\quad a\le 0.\\
    \end{aligned}
    \right.
        \end{equation*}
    Moreover, for $a>0$,
    \begin{flalign*}
         &\  W'_{a,v}(0,1,\Phi_v) &
    \end{flalign*}
    \begin{equation*}
         =
         \frac{\gamma_{2n+2}e^{-2\pi a}\cdot(2a)^n\pi^{n+1}}{\Gamma(1+n)}\cdot\Big(\sum_{i=1}^n\binom{n}{i}(4a\pi)^{-i}\Gamma(i)+
        \log (a\pi)+(\gamma-\sum_{i=1}^{n}\frac{1}{i})\Big).
    \end{equation*}
    Here $\gamma_{2n+2}=e^\frac{2\pi i(2n+2)}{8}$ is the Weil index.
        \item  The holomorphic projection
        \begin{equation*}
            \Pr'_\psi W'_{1,v}(0,g,\Phi_v)=c_3 W_{1,v}(0,g,\Phi_v),
        \end{equation*}
        where
        \begin{equation*}
            c_3=\sum_{i=1}^{n-1}\frac{n}{2i(n-i)}- \frac{1}{2}(\log4+\sum_{i=1}^{n}\frac{1}{i})=\sum_{i=1}^{n-1}\frac{1}{2i}-\frac{1}{2n}-\frac{\log4}{2}.
        \end{equation*}
    \end{enumerate}
\end{lemma}

\begin{proof}
    We again refer to \cite[Proposition 2.11]{YZZ2} for general result. Then we have the expression
    of $W_{a,v}(0,1,\Phi_v)$ directly. When $a>0$, note that from \cite[Proposition 2.11]{YZZ2}
    we have
    \begin{equation*}
    \begin{aligned}
        W_{a,v}(s,1,\Phi_v)=&\gamma_{2n+2}\frac{2\pi^{s+n+1}}{\Gamma(\frac{s}{2})\Gamma(\frac{s}{2}+n+1)}\int_{t>2a}e^{-2\pi (t-a)}t^{\frac{s}{2}+n}(t-2a)^{\frac{s}{2}-1}dt \\
        =&\gamma_{2n+2}\frac{2\pi^{s+n+1}}{\Gamma(\frac{s}{2})\Gamma(\frac{s}{2}+n+1)}e^{-2\pi a}\int_{t>0}e^{-2\pi t}(t+2a)^{\frac{s}{2}+n}t^{\frac{s}{2}-1}dt\\
        =&\gamma_{2n+2}e^{-2\pi a}\frac{\pi^{s+n+1}s}{\Gamma(\frac{s}{2}+1)\Gamma(\frac{s}{2}+n+1)}\int_{t>0}e^{-2\pi t}(t+2a)^{\frac{s}{2}+n}t^{\frac{s}{2}-1}dt .
    \end{aligned}
    \end{equation*}

    To see the derivative at $s=0$, note that the last integral is not convergent at $s=0$ due
    to the singularity of $t^{\frac{s}{2}-1}$ at $t=0$. The difference
    \begin{equation*}
        \begin{aligned}
            &\int_{t>0}e^{-2\pi t}(t+2a)^{\frac{s}{2}+n}t^{\frac{s}{2}-1}dt-\int_{t>0}e^{-2\pi t}(2a)^{\frac{s}{2}+n}t^{\frac{s}{2}-1}dt \\
            =&\int_{t>0}e^{-2\pi t}\frac{(t+2a)^{\frac{s}{2}+n}-(2a)^{\frac{s}{2}+n}}{t}t^{\frac{s}{2}}dt\\
            =&\int_{t>0}e^{-t}\Big(\sum_{i=1}^{n}\binom{n}{i}(2a)^{n-i}\cdot(\frac{t}{2\pi})^i\Big)\frac{dt}{t}+O(s),\\
            =&(2a)^n\sum_{i=1}^{n}\binom{n}{i}(4a\pi)^{-i}\cdot\Gamma(i)+O(s).
        \end{aligned}
    \end{equation*}

    Using this expression, we then have
    \begin{equation*}
    \begin{aligned}
        W_{a,v}(s,1,\Phi_v)=&\gamma_{2n+2}e^{-2\pi a}\frac{\pi^{s+n+1}s}{\Gamma(\frac{s}{2}+1)\Gamma(\frac{s}{2}+n+1)}\\
        &\cdot\Big(\sum_{i=0}^{n-1}\binom{n}{i}(2a)^i\cdot(2\pi)^{-n+i}\Gamma(n-i)+(2a)^{\frac{s}{2}+n}\cdot(2\pi)^{-\frac{s}{2}}\Gamma(\frac{s}{2})\Big)+O(s^2) \\
        =&\gamma_{2n+2}e^{-2\pi a}\frac{\pi^{s+n+1}}{\Gamma(\frac{s}{2}+1)\Gamma(\frac{s}{2}+n+1)} \\
        &\cdot\Big(\sum_{i=0}^{n-1}\binom{n}{i}(2a)^i\cdot(2\pi)^{-n+i}\Gamma(n-i)s+2(2a)^{\frac{s}{2}+n}\cdot(2\pi)^{-\frac{s}{2}}\Gamma(\frac{s+2}{2})\Big)+O(s^2)
    \end{aligned}
    \end{equation*}
    Taking $\lim_{s\rightarrow0}$ we recover $W_{a,v}(s,1,\Phi_v)$. Furthermore, taking the first derivative we have
    \begin{equation*}
    \begin{aligned}
        W'_{a,v}(0,1,\Phi_v)=\frac{\gamma_{2n+2}e^{-2\pi a}\cdot(2a)^n\pi^{n+1}}{\Gamma(1+n)}\cdot\Big(\sum_{i=1}^n\binom{n}{i}(4a\pi)^{-i}\Gamma(i)+
        \log (a\pi)-\psi^{(0)}(1+n)\Big).
    \end{aligned}
    \end{equation*}
    Here $\psi^{0}(z)$ is the polygamma function, i.e., $\psi^{0}(z)=\log(\Gamma(z))'$. It is well known that
    \begin{equation*}
        \psi^{0}(1+n)=\sum_{i=1}^{n}\frac{1}{i}-\gamma.
    \end{equation*}
    Thus we proved (1).

    Apply the proposition of Whittaker function in Section \ref{Theta}
    \begin{equation*}
    W_{1,v}(s,m(\sqrt{a}),\Phi_v)=a^{\frac{1-n}{2}-\frac{s}{2}}W_{a,v}(s,1,\Phi_v),
    \end{equation*}
    we then have
    \begin{equation*}
    \begin{aligned}
        &W_{1,v}'(0,m(\sqrt{a}),\Phi_v) \\
        &=\frac{\gamma_{2n+2}2^n\cdot\pi^{n+1}a^{\frac{n+1}{2}}\cdot e^{-2\pi a}}{\Gamma(1+n)}\Big(\sum_{i=1}^{n}\binom{n}{i}(4a\pi)^{-i}\Gamma(i)+ \log \pi-\psi^{(0)}(1+n)\Big).
    \end{aligned}
    \end{equation*}
    Note that this result agrees with the result in \cite[Lemma 3.3]{Yuan1}
    when $n=1$. Once again, the difference is that we have to use $m(\sqrt{a})$ to replace the element $d(a)\in \GL_2$ there, which makes our expression a little different. Nonetheless, we can still compute the holomorphic projection in the same way.

    To compute the holomorphic projection, we have
    \begin{equation*}
        \begin{aligned}
        &\Pr'_\psi W'_{1,v}(0,g,\Phi_v) \\
            &=(4\pi)^{n}\Gamma(n)^{-1}W^{(\mm_v)}(g)\cdot\widetilde{\lim_{s\rightarrow0}}\int_{Z(\RR)N(\RR)\backslash\UU(\RR)}\delta(h)^s W'_{1,v}(0,h,\Phi_v)
        \overline{W^{(\mm_v)}(h)}dh\\
        &=(4\pi)^{n}\Gamma(n)^{-1}W^{(\mm_v)}(g)\cdot\widetilde{\lim_{s\rightarrow0}}\int_{t>0}t^{s+\frac{n-1}{2}}e^{-2\pi t}W'_{1,v}(0,m(\sqrt{t}),\Phi_v)\frac{dt}{t} \\
        &=\gamma_{2n+2}(4\pi)^{n}\Gamma(n)^{-1}W^{(\mm_v)}(g)\cdot2^n\pi^{n+1}\cdot\Gamma(1+n)^{-1} \\
        &\cdot\widetilde{\lim_{s\rightarrow0}}\int_{t>0}t^{s+n}e^{-4\pi t}\Big(\sum_{i=1}^{n}\binom{n}{i}(4t\pi)^{-i}\Gamma(i)+ \log \pi-\psi^{(0)}(1+n)\Big)\frac{dt}{t}\\
        &=\gamma_{2n+2}\Gamma(n)^{-1}W^{(\mm_v)}(g)\cdot2^n\pi^{n+1}\cdot\Gamma(1+n)^{-1} \\
        &\cdot\widetilde{\lim_{s\rightarrow0}}(4\pi)^{-s}\Big(\sum_{i=1}^{n}\binom{n}{i}\Gamma(n-i+s)\Gamma(i)+(\log \pi-\psi^{(0)}(1+n))\cdot\Gamma(n+s)\Big)
        \end{aligned}
    \end{equation*}
    Here we use the Iwasawa decomposition at the second equation. Taking care
    of the $(4\pi)^{-s}\cdot\Gamma(s)$-term, we know the
    constant term of the above sum in limit is equal to
    \begin{equation*}
        \begin{aligned}
        \sum_{i=1}^{n-1}\binom{n}{i}\Gamma(n-i)\Gamma(i)-\Gamma(n)\cdot\big(\log4+\gamma+\psi^{(0)}(1+n)\big).
        \end{aligned}
    \end{equation*}

    Moreover, it is not hard to compute from (1) that
    \begin{equation*}
        W_{1,v}(0,g,\Phi_v)=\gamma_{2n+2}(2\pi)^{1+n}\Gamma(1+n)^{-1}W^{(\mm_v)}(g).
    \end{equation*}
    Thus, we conclude that
    \begin{equation*}
        \Pr'_\psi W'_{1,v}(0,g,\Phi_v)=
        \Big(\sum_{i=1}^{n-1}\frac{n}{2i(n-i)}- \frac{1}{2}(\log4+\sum_{i=1}^{n}\frac{1}{i}) \Big)W_{1,v}(0,g,\Phi_v).   \end{equation*}

    Note that when $n=1$, our result is compatible with \cite[Lemma 3.3]{Yuan1}.
\end{proof}

As we mentioned in the introduction, this constant $c_3$ is exactly the constant $\mathfrak{b}$ (up to a factor of 2) in \cite[Thm 1.1.1]{Qiu} that was not made explicit.

\subsection{Explicit computations of J-part: non-archimedean places}
Denote by $p_v$ the maximal ideal of $\mathcal{O}_{\Fn}$, and denotes by
$d_v\in F_v$ the local different of $F_v/\QQ_p$. We also denote by $D_v$ the relative discriminant $D_{E_v/F_v}$. Following the assumption
of Schwartz function $\Phi_v$ in Section \ref{Choice of Schwartz function}, we have the following explicit results of the local Whittaker function associated with $\Phi_v$.

\begin{lemma}\label{Explicit nonarchimedean}
    Let $v$ be a non-archimedean place of $F$, and let $a\in\Fn^\times$. Denote by  $r=v(a)$.
    \begin{enumerate}
        \item If $v$ is unramified in $E/F$, then $W_{a,v}(s,1,\Phi_v)$ is nonzero only if $v(a)\ge -v(d_v)$. In this case, for $a\in\mathcal{O}_{\Fn}$, when $2\nmid n$ or $v$ is split in $E/F$, we have
        \begin{equation*}
            \begin{aligned}
                W_{a,v}(s,1,\Phi_v)=&|d_v|^{s+n+\frac{1}{2}}\frac{(1-N_v^{-(s+n+1)})(1-N_v^{-(r+1)(s+n)})}{1-N_v^{-(s+n)}}\\
                &+|d_v|^{n+\frac{3}{2}}\frac{(1-N_v^{-s})(1-|d_v|^{s-1})}{1-N_v^{-(s-1)}},
            \end{aligned}
        \end{equation*}
        and when $2|n$ and $v$ is inert in $E/F$, we have
        \begin{equation*}
            \begin{aligned}
                W_{a,v}(s,1,\Phi_v)=&|d_v|^{s+n+\frac{1}{2}}\frac{(1+N_v^{-(s+n+1)})(1-(-1)^{r+1}N_v^{-(r+1)(s+n)})}{1+N_v^{-(s+n)}}\\
                &+|d_v|^{n+\frac{3}{2}}\frac{(1-N_v^{-s})(1-|d_v|^{s-1})}{1-N_v^{-(s-1)}}.
            \end{aligned}
        \end{equation*}

        Therefore, for $a\in\mathcal{O}_{\Fn}$, when $2\nmid n$ or $v$ is split in $E/F$,
        \begin{equation*}
            \begin{aligned}
                &W'_{a,v}(0,1,\Phi_v)-\frac{1}{2}\log|a|_v W_{a,v}(0,1,\Phi_v)\\
                =&\big(-\frac{L'(n+1,\eta_v^{n+1})}{L(n+1,\eta_v^{n+1})}+\log|d_v|\big)W_{a,v}(0,1,\Phi_v)\\
            &+|d_v|^{n+\frac{1}{2}}\frac{1-N_v^{-(n+1)}}{2(1-N_v^{-n})^2} \big(r(1-N_v^{-(r+2)n})-(r+2)(N_v^{-n}-N_v^{-(r+1)n})\big)\log N_v\\
            &+|d_v|^{n+\frac{1}{2}}\frac{1-|d_v|}{N_v-1}\log N_v,
            \end{aligned}
         \end{equation*}
        and when $2|n$ and $v$ is inert in $E/F$,
        \begin{equation*}
            \begin{aligned}
                &W'_{a,v}(0,1,\Phi_v)-\frac{1}{2}\log|a|_v W_{a,v}(0,1,\Phi_v)\\
                =&\big(-\frac{L'(n+1,\eta_v^{n+1})}{L(n+1,\eta_v^{n+1})}+\log|d_v|\big)W_{a,v}(0,1,\Phi_v)\\
            +&|d_v|^{n+\frac{1}{2}}\frac{1+N_v^{-(n+1)}}{2(1+N_v^{-n})^2} \big(r(1-(-1)^{r}N_v^{-(r+2)n})+(r+2)(N_v^{-n}-(-1)^{r}N_v^{-(r+1)n})\big)\log N_v\\
            +&|d_v|^{n+\frac{1}{2}}\frac{1-|d_v|}{N_v-1}\log N_v.
            \end{aligned}
        \end{equation*}
        Note that
        \begin{equation*}
            \frac{\zeta'_v(n+1)}{\zeta_v(n+1)}=\frac{L'(n+1,\eta_v^{n+1})}{L(n+1,\eta_v^{n+1})}
        \end{equation*}
        when $n$ is odd or $v$ is split in $E/F$.
        \item If $v$ is ramified in $E/F$, then $W_{a,v}(s,1,\Phi_v)$ is nonzero only if $v(a)\ge 0$. In this case, for $a\in\mathcal{O}_{\Fn}$, if $2\nmid n$, we have
        \begin{equation*}
        \begin{aligned}
            &W'_{a,v}(0,1,\Phi_v)-\frac{1}{2}\log|a|_v W_{a,v}(0,1,\Phi_v)\\
            =&-\frac{\zeta'_v(n+1)}{\zeta_v(n+1)}W_{a,v}(0,1,\Phi_v)-\frac{r}{2}(1-N_v^{-(n+1)})\log N_v+N_v^{-(n+1)}\log N_v\\
            &+\frac{1-N_v^{-(n+1)}}{2(1-N_v^{-n})^2}\big(r-(r+2)N_v^{-n}+(r+2)N_v^{-(r+1)n}-rN_v^{-(r+2)n}\big)\log N_v,
        \end{aligned}
    \end{equation*}
        If $2|n$ and $a\notin \Nm(E_v^\times)\cdot(-1)^\frac{n}{2}d_\BV$, we have
        \begin{equation*}
            \begin{aligned}
                W'_{a,v}(0,1,\Phi_v)-\frac{1}{2}\log|a|_v W_{a,v}(0,1,\Phi_v)
                =\big((\frac{r+2}{2})N_v^{-(r+1)n-\frac{1}{2}}+\frac{r}{2}N_v^{-(n+\frac{1}{2})}\big)\log N_v.
            \end{aligned}
        \end{equation*}
        If $2|n$ and $a\in \Nm(E_v^\times)\cdot(-1)^\frac{n}{2}d_\BV$, we have
        \begin{equation*}
            \begin{aligned}
                W'_{a,v}(0,1,\Phi_v)-\frac{1}{2}\log|a|_v W_{a,v}(0,1,\Phi_v)
                =\big(-(\frac{r+2}{2})N_v^{-(r+1)n-\frac{1}{2}}+\frac{r}{2}N_v^{-(n+\frac{1}{2})}\big)\log N_v.
            \end{aligned}
        \end{equation*}
        Note that when $2|n$, $L(n+1,\eta_v^{n+1})=1$, and for those $a$ with $r=0$, $a\in q(\Lambda_v)$ if and only if $a\in \Nm(E_v^\times)\cdot(-1)^\frac{n}{2}d_\BV$.
    \end{enumerate}
\end{lemma}

\begin{proof}
    The proof follows the proof in \cite[Lemma 3.4]{Yuan1}. We first deal with the case of $d_v=1$, then we consider a general formula for Whittaker function. Finally, we use this general formula to compute the case of $d_v\ne 1$. In fact, if $|d_v|\ne1$, $\Phi_v$ is not invariant under the Weil representation, which causes some difficulty.

    \textbf{Step 1. Unramified case:} In this case, $v$ is unramified in both $E/F$ and $F/\QQ$ with $\epsilon(\BV_v)=1$.
    We recall again the definition of our local Whittaker function
    \begin{equation*}
        W_{a,v}(s,g,\Phi_v)=\int_{F_v}\delta_v(wn(b)g)^s r(wn(b)g)\Phi_v(0)\psi(-ab)db.
    \end{equation*}
    Note that in our case, $\Phi_v$ is invariant under the action of $r(m(a))$ for any $a\in\mathcal{O}_{E_v}^\times$ and $r(wn(b))$ for any $b\in\mathcal{O}_{F_v}$ from the unramified assumption. Especially, note that the Weil index is always trivial in this case. Also recall that the Iwasawa decomposition gives an equation
    \begin{equation}\label{matrix equation}
        wn(b)=n(-b^{-1})m(b^{-1})wn(-b^{-1})w.
    \end{equation}
    This equation then implies
    \begin{equation*}
        r(wn(b))\Phi_v(0)=\chi_{\BV,v}(b^{-1})\delta(wn(b))^{n+1},\quad b^{-1}\in \mathcal{O}_{F_v},
    \end{equation*}
    Here $\chi_{\BV,v}|_{F_v^\times}=\eta_v^{n+1}$. Thus it is trivial at a split place and also at an inert place when $n$ is odd, whereas at an inert place with $n$ even one has $\chi_{\BV,v}(b^{-1})=(-1)^{v(b)}$. This is the character denoted by $\chi_v$ in the raw local calculation below.

    As a consequence, we have
    \begin{equation*}
    \begin{aligned}
        W_{a,v}(s,1,\Phi_v)=\int_{\mathcal{O}_{\Fn}}\psi(-ab)db+\int_{\Fn-\mathcal{O}_{\Fn}}
        \chi_{\BV,v}(b^{-1})|b|^{-(s+n+1)}\psi(-ab)db.
    \end{aligned}
    \end{equation*}
    Write the domain $\Fn-\mathcal{O}_{\Fn}$ as the disjoint union of $p_v^{-k}-p_v^{-(k-1)}$ for $k\ge 1$, recall the Haar measure we chosen in Section \ref{Notation}, for $a\in\mathcal{O}_{\Fn}$, we conclude that if $2\nmid n$ or $v$ is split,
    \begin{equation*}
        W_{a,v}(s,1,\Phi_v)=(1-N_v^{-(s+n+1)})\sum_{m=0}^{v(a)} N_v^{-m(s+n+1)}\cdot N_v^m=\frac{(1-N_v^{-(s+n+1)})(1-N_v^{-(v(a)+1)(s+n)})}{1-N_v^{-(s+n)}},
    \end{equation*}
    while if $2|n$ and $v$ is inert,
    \begin{equation*}
        W_{a,v}(s,1,\Phi_v)=\frac{(1+N_v^{-(s+n+1)})(1-(-1)^{v(a)+1}N_v^{-(v(a)+1)(s+n)})}{1+N_v^{-(s+n)}}.
    \end{equation*}
    For other $a$, $W_{a,v}(s,1,\Phi_v)=0$.

    It follows that
    \begin{equation*}
        W_{a,v}(0,1,\Phi_v)=\frac{(1-N_v^{-(n+1)})(1-N_v^{-(v(a)+1)n})}{1-N_v^{-n}}
    \end{equation*}
    when $2\nmid n$ or $v$ is split, and
    \begin{equation*}
        W_{a,v}(0,1,\Phi_v)=\frac{(1+N_v^{-(n+1)})(1-(-N_v^{-n})^{(v(a)+1)})}{1+N_v^{-n}}
    \end{equation*}
    otherwise. Moreover, a direct computation shows that if $2\nmid n$ or $v$ is split
    \begin{equation*}
        \begin{aligned}
            &W'_{a,v}(0,1,\Phi_v)-\frac{1}{2}\log|a|_v W_{a,v}(0,1,\Phi_v)\\
            =&W_{a,v}(0,1,\Phi_v)\big(\frac{N_v^{-(n+1)}}{1-N_v^{-(n+1)}}+\frac{(v(a)+1)N_v^{-(v(a)+1)n}}{1-N_v^{-(v(a)+1)n}}-\frac{N_v^{-n}}{1-N_v^{-n}}+\frac{1}{2}v(a)\big)\log N_v \\
            =&-\frac{\zeta'_v(n+1)}{\zeta_v(n+1)}W_{a,v}(0,1,\Phi_v)\\
            &+\frac{1-N_v^{-(n+1)}}{2(1-N_v^{-n})^2}\big(r-(r+2)N_v^{-n}+(r+2)N_v^{-(r+1)n}-rN_v^{-(r+2)n}\big)\log N_v,
        \end{aligned}
    \end{equation*}
    if $2|n$ and $v$ is inert
    \begin{equation*}
        \begin{aligned}
            &W'_{a,v}(0,1,\Phi_v)-\frac{1}{2}\log|a|_v W_{a,v}(0,1,\Phi_v)\\
            =&-\frac{L'(n+1,\eta_v)}{L(n+1,\eta_v)}W_{a,v}(0,1,\Phi_v)\\
            &+\frac{1+N_v^{-(n+1)}}{2(1+N_v^{-n})^2}\big(r+(r+2)N_v^{-n}+(-1)^{r+1}(r+2)N_v^{-(r+1)n}+(-1)^{r+1}rN_v^{-(r+2)n}\big)\log N_v,
        \end{aligned}
    \end{equation*}
    Here recall that $\eta$ is the associated quadratic Hecke character of $F^\times\backslash\BA_F^\times$ via the class field theory, and we use the result
    \begin{equation*}
        \frac{\zeta'_v(s)}{\zeta_v(s)}=\frac{-\log N_v\cdot N_v^{-s}}{1-N_v^{-s}},\quad \frac{L'(s,\eta_v)}{L(s,\eta_v)}=\frac{\log N_v \cdot N_v^{-s}}{1+N_v^{-s}}
    \end{equation*}
    This proves part (1) when $|d_v|=1$.

    \textbf{Step 2. A general formula:} Following the result and proof in
    \cite[Proposition 6.10(1)]{YZZ2}, there is a general formula holds for
    non-trivial $d_v$:
    \begin{equation*}
        W_{a,v}(s,1,\Phi_v)=\gamma(\BV_v,q)|d_v|^\frac{1}{2}(1-N_v^{-s})\sum_{m=0}^\infty N_v^{-m(s-1)}\int_{D_m(a)}\Phi_v(x)\de x.
    \end{equation*}
    Here $\de x$ is the self-dual measure of $(\BV_v,q)$, and
    \begin{equation*}
        D_m(a)=\{x\in\BV_v|q(x)\in a+p_v^m d_v^{-1}\}
    \end{equation*}
    is a subset of $\BV_v$. The Weil index comes from the action of $w$. There is no essential difference if we replace the underlying quadratic space by Hermitian space, and all the notations are compatible.

    \textbf{Step 3. Inert case:} Now we assume that $v$ is inert in $E/F$, and $d_v$ is arbitrary in Step 1. The major difficulty is that, for $|d_v|\ne 1$, the characteristic function $1_{\mathcal{O}_{\Fn}}$ is not self-dual under the normalized character $\psi_v: F_v\rightarrow\CC^\times$ mentioned in Section \ref{Notation}. Consequently, the Schwartz function $\Phi_v$ is not invariant under the action of $\UU(\mathcal{O}_{\Fn})$, and the method in Step 1 does not work any more. Thus, we are going to use the formula in Step 2.

    We have
    \begin{equation*}
         W_{a,v}(s,1,\Phi_v)=|d_v|^\frac{1}{2}(1-N_v^{-s})\sum_{m=0}^\infty N_v^{-m(s-1)}\vol(D_m(a)\cap\Lambda_v).
    \end{equation*}
    Note that in this case, $\vol(\Lambda_v)=|d_v|^{n+1}$ following the volume of $\BV_v$ chosen in Section \ref{Notation}.  We then write
    \begin{equation*}
        W_{a,v}(s,1,\Phi_v)=|d_v|^{n+\frac{3}{2}}(1-N_v^{-s})\sum_{m=0}^\infty \frac{N_v^{-m(s-1)}\vol(D_m(a)\cap\Lambda_v)}{\vol(\Lambda_v)}.
    \end{equation*}
    Split this summation according to $m<v(d_v)$ and $m\ge v(d_v)$, since we assume $a\in \mathcal{O}_{\Fn}$, it follows that
    \begin{equation*}
    \begin{aligned}
        W_{a,v}(s,1,\Phi_v)_{m< v(d_v)}=&|d_v|^{n+\frac{3}{2}}(1-N_v^{-s})\sum_{m=0}^{v(d_v)-1} N_v^{-m(s-1)}\\
        =&|d_v|^{n+\frac{3}{2}}\frac{(1-N_v^{-s})(1-|d_v|^{s-1})}{(1-N_v^{-(s-1)})}.
    \end{aligned}
    \end{equation*}
    This is because $\Lambda_v\subset D_m(a)$ in this case.

    Meanwhile, a direct calculation of $W_{a,v}(s,1,\Phi_v)_{m\ge v(d_v)}$ is quite involved, so we compare it with the unramified case instead. Assume $2\nmid n$ here,  denote by
    \begin{equation*}
        D_m(a)^\circ=\{x\in\BV_v|q(x)\in a+p_v^m\},
    \end{equation*}
    which is equal to the set $D_m(a)$ in the unramified case in Step 1. For $m\ge v(d_v)$, the substitution $m\mapsto m+v(d_v)$ gives
    \begin{equation*}
        W_{a,v}(s,1,\Phi_v)_{m\ge v(d_v)}=|d_v|^{s+n+\frac{1}{2}}(1-N_v^{-s})\sum_{m=0}^\infty \frac{N_v^{-m(s-1)}\vol(D_m(a)^\circ\cap\Lambda_v)}{\vol(\Lambda_v)}.
    \end{equation*}
    This is equal to $W_{a,v}(s,1,\Phi_v)$ in the case $|d_v|=1$ considered in Step 1. In other words, the result in Step 1 implies
    \begin{equation*}
        (1-N_v^{-s})\sum_{m=0}^\infty \frac{N_v^{-m(s-1)}\vol(D_m(a)^\circ\cap\Lambda_v)}{\vol(\Lambda_v)}=\frac{(1-N_v^{-(s+n+1)})(1-N_v^{-(v(a)+1)(s+n)})}{1-N_v^{-(s+n)}}.
    \end{equation*}
    Hence, we have
    \begin{equation*}
        W_{a,v}(s,1,\Phi_v)_{m\ge v(d_v)}=|d_v|^{s+n+\frac{1}{2}}\frac{(1-N_v^{-(s+n+1)})(1-N_v^{-(v(a)+1)(s+n)})}{1-N_v^{-(s+n)}}.
    \end{equation*}
    Thus, we can combine these two parts to get a formula of $W_{a,v}(s,1,\Phi_v)$, and it is routine to check part (1) when $2\nmid n$. When $2|n$, there is no essential difference, we have
    \begin{equation*}
        W_{a,v}(s,1,\Phi_v)_{m\ge v(d_v)}=|d_v|^{s+n+\frac{1}{2}}\frac{(1+N_v^{-(s+n+1)})(1-(-1)^{v(a)+1}N_v^{-(v(a)+1)(s+n)})}{1+N_v^{-(s+n)}}.
    \end{equation*}
    by similar discussion. Hence we complete all the proof of part (1).

    \textbf{Step 4. Ramified case:} Now assume that $v$ is ramified in $E/F$. By the standing hypothesis of Section \ref{Choice of Schwartz function}, $|d_v|=1$. Note that by our choice the Hermitian lattice $\Lambda_v$ is $\varpi_{E_v}$-modular or almost $\varpi_{E_v}$-modular, and hence $\Phi_v$ is not invariant under the action of $\UU(\mathcal{O}_{F_v})$. Thus, the method used in Step 1 cannot be directly applied in this context. Meanwhile, if $2\nmid n$, unlike our discussion in Step 3, we do not need to use the general formula in Step 2. In fact, we can check the Weil representation carefully in this case, and then modify our proof in Step 1 to get the correct proof for $2\nmid n$ case. While for $2|n$ case, we have to use the general formula in Step 2, and reduce to $2\nmid n$ case by some elementary computation.

    Recall the definition of local Whittaker function in Step 1. Although $\Lambda_v$ is no longer self-dual in this case, the Schwartz function $\Phi_v$ is still invariant under the action of $r(m(a))$ and $r(n(b))$ for any $a\in\mathcal{O}_{E_v}^\times$ and $b\in \mathcal{O}_{F_v}$. Furthermore, we have
    \begin{equation*}
        r(w)\Phi_v=\vol(\Lambda_v)\Phi'_v,
    \end{equation*}
    where $\Phi_v'$ is the characteristic function of lattice $\varpi_{E_v}^{-1}\Lambda_v^\vee$. Here we use the fact that $\gamma(\BV_v,q)=1$ in ramified case and
    \begin{equation*}
        \varpi_{E_v}^{-1}\Lambda_v^\vee=\{x\in \BV_v\Big|\tr_{E_v/F_v}(\langle x,y\rangle)\in\mathcal{O}_{F_v},\ \forall y\in\Lambda_v\}.
    \end{equation*}

    For further discussion, we need to consider the cases separately based on the parity of $n$. First assume $2\nmid n$, in this case $\Lambda_v^\vee=\varpi_{E_v}^{-1}\Lambda_v$, and $q(\Lambda_v)=p_v\mathcal{O}_{F_v}$. Following the discussion in Step 1, on the one hand, it is not hard to check that if $v(b)\ge 0$
    \begin{equation*}
        r(wn(b))\Phi_v(0)=r(w)\Phi_v(0)=\vol(\Lambda_v)=N_v^{-(n+1)}.
    \end{equation*}
    On the other hand, if $v(b)<0$, apply the matrix identity \eqref{matrix equation}, the relation $r(w)\Phi_v=\vol(\Lambda_v)\Phi_v'$, and the invariance of $\Phi_v'$ under the relevant integral unipotent factor. Evaluating at the origin and using that $\chi_{\BV,v}|_{F_v^\times}=1$ in the present parity gives
\begin{equation*}
 r(wn(b))\Phi_v(0)=\delta(wn(b))^{n+1}=|b|^{-(n+1)}.
\end{equation*}

    Thus, combine these two parts, we have
    \begin{equation*}
    \begin{aligned}
        W_{a,v}(s,1,\Phi_v)=\vol(\Lambda_v)\cdot\int_{\mathcal{O}_{\Fn}}\psi(-ab)db
        +\int_{F_v-\mathcal{O}_{\Fn}}|b|^{-(s+n+1)}\psi(-ab)db.
    \end{aligned}
    \end{equation*}
    Then the remaining computation is the same as the computation in Step 1, since the second term here is exactly the same with the corresponding term in unramified case. As a conclusion, for $a\in\mathcal{O}_{F_v}$  and $r=v(a)$, we have
    \begin{equation*}
       W_{a,v}(s,1,\Phi_v)=\frac{(1-N_v^{-(s+n+1)})(1-N_v^{-(r+1)(s+n)})}{1-N_v^{-(s+n)}}-(1-N_v^{-(n+1)}).
    \end{equation*}
    Moreover, we have
    \begin{equation*}
        \begin{aligned}
            &W'_{a,v}(0,1,\Phi_v)-\frac{1}{2}\log|a|_v W_{a,v}(0,1,\Phi_v)\\
            =&-\frac{\zeta'_v(n+1)}{\zeta_v(n+1)}W_{a,v}(0,1,\Phi_v)-\frac{r}{2}(1-N_v^{-(n+1)})\log N_v+N_v^{-(n+1)}\log N_v\\
            &+\frac{1-N_v^{-(n+1)}}{2(1-N_v^{-n})^2}\big(r-(r+2)N_v^{-n}+(r+2)N_v^{-(r+1)n}-rN_v^{-(r+2)n}\big)\log N_v,
        \end{aligned}
    \end{equation*}
    which can be easily derived from the formula of unramified case.

    It remains to consider the case when $2|n$. In this case $\vol(\Lambda_v)=N_v^{-(n+\frac{1}{2})}$. By the local classification, choose an orthogonal decomposition $\Lambda_v=\Lambda_v^{(n)}\oplus\Lambda_v^{(1)}$, where $\Lambda_v^{(n)}$ is $\varpi_{E_v}$-modular of rank $n$ and $\Lambda_v^{(1)}$ is a self-dual Hermitian line. This is a chosen standard decomposition, not a canonical one. For convenience, we also denote by $\Phi_v^{(n)}$ the characteristic function of $\Lambda_v^{(n)}$. Under this decomposition, we have
    \begin{equation*}
\varpi_{E_v}^{-1}\Lambda_v^\vee=\varpi_{E_v}^{-2}\Lambda_v^{(n)}\oplus\varpi_{E_v}^{-1}\Lambda_v^{(1)}.
    \end{equation*}
    This implies that if $v(b)<0$, $\Phi_v'$ is still invariant under the action of $r(n(-b^{-1}))$. Moreover, in the present case $n$ is even, so $\chi_{\BV,v}(b^{-1})=\eta_v(b^{-1})$; this equals $1$ exactly when $b\in\Nm(E_v^\times)$ and equals $-1$ otherwise. This will lead to significant computational challenges because the value of the character will no longer solely depend on the valuation. Thus, we instead turn to employing the method described in Step 2, i.e., we consider the formula
     \begin{equation*}
        W_{a,v}(s,1,\Phi_v)=\vol(\Lambda_v)(1-N_v^{-s})\sum_{m=0}^\infty \frac{N_v^{-m(s-1)}\vol(D_m(a)\cap\Lambda_v)}{\vol(\Lambda_v)}.
    \end{equation*}

    By definition, we have
    \begin{equation*}
        q(\Lambda_v)=p_v\mathcal{O}_{F_v}+q(\Lambda^{(1)}_v).
    \end{equation*}
    Thus, on the one hand, if $v(a)=0$ and $a\notin q(\Lambda^{(1)}_v)$, we have
    \begin{equation*}
        D_m(a)\cap \Lambda_v=\left\{
    \begin{aligned}
        \nonumber
        &\Lambda_v \ \ \ m=0;\\
        &\emptyset\ \ \ m\ge 1.\\
    \end{aligned}
    \right.
    \end{equation*}
    This implies
    \begin{equation*}
        W_{a,v}(s,1,\Phi_v)=N_v^{-(n+\frac{1}{2})}-N_v^{-(s+n+\frac{1}{2})}.
    \end{equation*}
    Note that by computing the Hermitian determinant, $q(\Lambda_v^{(1)})=\Nm(\mathcal{O}_{E_v})\cdot (-1)^\frac{n}{2}d_\BV$.

    On the other hand, if $v(a)=0$ and $a\in q(\Lambda_v^{(1)})$, the situation is highly delicate. We claim that there is an equation
    \begin{equation*}
        W_{a,v}(s,1,\Phi_v)=\sum_{m=0}^\infty C^{0}_m\frac{\vol(\Lambda_v)}{\vol(\Lambda_v^{(n)})}W_{p_v^{m},v}(s,1,\Phi_v^{(n)}),
    \end{equation*}
    where
    \begin{equation*}
        C^{0}_m=\left\{
    \begin{aligned}
        \nonumber
        &1-2N_v^{-1} \ \ \ m=0;\\
        &2N_v^{-m}-2N_v^{-(m+1)} \ \ \ m\ge 1,\\
    \end{aligned}
    \right.
    \end{equation*}
    and $W_{p_v^{m},v}(s,1,\Phi_v^{(n)})$ is the local Whittaker function associated with $\Phi_v^{(n)}$.

    Indeed, it is obvious to see that for such $a$,
    \begin{equation*}
        C^{0}_m=\frac{\vol(\{x\in\Lambda_v^{(1)}\big|\ v(a-q(x))=m\})}{\vol(\Lambda_v^{(1)})}.
    \end{equation*}
    This is also the reason why we use the notation $C_m^{0}$, where the superscript $0$ is the valuation of $a$. Then the computation of $W_{a,v}(s,1,\Phi_v)$ in this case can be reduced to the computation of each $W_{a-q(x),v}(s,1,\Phi_v^{(n)})$ for $x\in \Lambda_v^{(1)}$. This implies
    \begin{equation*}
        W_{a,v}(s,1,\Phi_v)=N_v^{-(n+\frac{1}{2})}+N_v^{-(s+n+\frac{1}{2})}.
    \end{equation*}

    It remains to consider the general case when $r=v(a)\ge 1$. Similar to $v(a)=0$ case, the result of $W_{a,v}(s,1,\Phi_v)$ depends on whether $a\in q(\Lambda_v^{(1)})$. In fact, we can generalize the above approach. On the one hand, if $a\notin q(\Lambda_v^{(1)})$, we have an equation
    \begin{equation}\label{computing Whittaker function using induction}
        W_{a,v}(s,1,\Phi_v)=\sum_{m=0}^r C^r_m\frac{\vol(\Lambda_v)}{\vol(\Lambda_v^{(n)})}W_{p_v^{m},v}(s,1,\Phi_v^{(n)}),
    \end{equation}
    where in this case
    \begin{equation*}
        C^r_m=\left\{
    \begin{aligned}
        \nonumber
        &N_v^{-m}-N_v^{-(m+1)} \ \ \ 0\le m<r;\\
        &N_v^{-r} \ \ \ m=r.\\
    \end{aligned}
    \right.
    \end{equation*}
    On the other hand, if $a\in q(\Lambda_v^{(1)})$, a similar induction formula as \ref{computing Whittaker function using induction} is written as
    \begin{equation*}
        W_{a,v}(s,1,\Phi_v)=\sum_{m=0}^\infty C^r_m\frac{\vol(\Lambda_v)}{\vol(\Lambda_v^{(n)})}W_{p_v^{m},v}(s,1,\Phi_v^{(n)}),
    \end{equation*}
    where in this case, we instead have
    \begin{equation*}
        C^r_m=\left\{
    \begin{aligned}
        \nonumber
        &N_v^{-m}-N_v^{-(m+1)} \ \ \ 0\le m<r;\\
        &N_v^{-r}-2N_v^{-(r+1)} \ \ \ m=r;\\
        &2N_v^{-m}-2N_v^{-(m+1)} \ \ \ m>r.
    \end{aligned}
    \right.
    \end{equation*}
    The proof methods for these two equations are consistent with the previous case of $r=0$.

    Hence, we conclude that in general, when $2|n$ and $r\ge 0$,
    \begin{equation*}
        W_{a,v}(s,1,\Phi_v)=\left\{
    \begin{aligned}
        \nonumber
        &N_v^{-(n+\frac{1}{2})}-N_v^{-(r+1)(s+n)-\frac{1}{2}} \ \ \ a\notin q(\Lambda_v^{(1)});\\
        &N_v^{-(n+\frac{1}{2})}+N_v^{-(r+1)(s+n)-\frac{1}{2}} \ \ \ a\in q(\Lambda_v^{(1)}).\\
    \end{aligned}
    \right.
    \end{equation*}
    The formula of $W'_{a,v}(0,1,\Phi_v)-\frac{1}{2}\log|a|_v W_{a,v}(0,1,\Phi_v)$ can be computed easily. Thus, we finish all the proof.

\end{proof}

From our proof process, it is evident that although our methods of computation follows from \cite[Lemma 3.4]{Yuan1}, there are some differences in both the techniques and the results. The major difference comes from our choice of the Hermitian lattice $\Lambda_v$ at those ramified prime, where the computation depends on the parity of $n$. When $n$ is odd, the computation is simpler, since the $\varpi_{E_v}$-modular lattice is self-dual after taking the trace map. See Remark \ref{dual by trace remark}; when $n$ is even, the expression of $W_{a,v}(s,1,\Phi_v)$ no longer solely relies on the valuation of $a$, which makes the computation more complicated.

It is worth noting that, in the specific calculation results mentioned above, terms such as
\begin{equation*}
    |d_v|^{n+\frac{1}{2}}\frac{1-N_v^{-(n+1)}}{2(1-N_v^{-n})^2} \big(r(1-N_v^{-(r+2)n})-(r+2)(N_v^{-n}-N_v^{-(r+1)n})\big)\log N_v
\end{equation*}
will magically reappear in the computation of height series in \cite{Guo1,Guo2}.

\begin{remark}\label{induction of Whittaker function remark}
    Note that the induction method used in Step 4 is valid for general cases. For instance, if $v$ is split in $E/F$, we can also use the formula \ref{computing Whittaker function using induction} to compute the expression of $W_{a,v}(s,1,\Phi_v)$. In this case, the coefficients $C_m^r$ are different. For example, it is not hard to check that
    \begin{equation*}
        C^0_m=\left\{
    \begin{aligned}
        \nonumber
        &1-N_v^{-1}+N_v^{-2} \ \ \ m=0;\\
        &N_v^{-m}-2N_v^{-(m+1)}+N_v^{-(m+2)} \ \ \ m>0.
    \end{aligned}
    \right.
    \end{equation*}
    Similarly, we can also use such induction method to compute the case when $v$ is inert in $E/F$.
\end{remark}

The following corollary provides the explicit expression of the corresponding Whittaker function at a ramified place when we choose different Schwartz functions. While these conclusions may not be used in the subsequent discussions, they will provide us with a better understanding of the properties of the Whittaker function.

\begin{corollary}\label{dual Schwartz function corollary}
    Suppose $v$ is ramified in $E/F$.
    \begin{enumerate}
        \item Suppose $2\nmid n$, and define $\Phi_v^\vee=1_{\Lambda_v^\vee}$.
              In this case, for $a\in\mathcal{O}_{F_v}$, we have
        \begin{equation*}
            W_{a,v}(s,1,\Phi_v^\vee)=\frac{(1-N_v^{-(s+n+1)})(1-N_v^{-(v(a)+1)(s+n)})}{1-N_v^{-(s+n)}},
        \end{equation*}
        and the difference
        \begin{equation*}
            W_{a,v}(s,1,\Phi_v^\vee)- W_{a,v}(s,1,\Phi_v)=1-N_v^{-(n+1)}
        \end{equation*}
        is a constant.
        \item Suppose $2\nmid n$ and the local Hermitian space $\BV_v$ admits a self-dual lattice; let $\Phi'_v$ be its characteristic function. In this case, for $a\in\mathcal{O}_{F_v}$, we have
        \begin{equation*}
            W_{a,v}(s,1,\Phi'_v)=\frac{(1-N_v^{-(s+n+1)})(1-N_v^{-(v(a)+1)(s+n)})}{1-N_v^{-(s+n)}}-(1-N_v^{-(\frac{n+1}{2})}),
        \end{equation*}
        and the difference
        \begin{equation*}
            W_{a,v}(s,1,\Phi'_v)- W_{a,v}(s,1,\Phi_v)=N_v^{-(\frac{n+1}{2})}-N_v^{-(n+1)}
        \end{equation*}
        is a constant.
        \item Suppose $2|n$, and $\Phi_v^\vee=1_{\Lambda_v^\vee}$. In this case, for $a\in\mathcal{O}_{F_v}$, we have
        \begin{equation*}
        W_{a,v}(s,1,\Phi^\vee_v)=\left\{
        \begin{aligned}
           \nonumber
           &N_v^{-\frac{1}{2}}-N_v^{-(r+1)(s+n)-\frac{1}{2}} \ \ \ a\notin q(\Lambda_v^{(1)});\\
           &N_v^{-\frac{1}{2}}+N_v^{-(r+1)(s+n)-\frac{1}{2}} \ \ \ a\in q(\Lambda_v^{(1)}),\\
        \end{aligned}
          \right.
       \end{equation*}
       and the difference
       \begin{equation*}
           W_{a,v}(s,1,\Phi_v^\vee)- W_{a,v}(s,1,\Phi_v)=N_v^{-\frac{1}{2}}-N_v^{-(n+\frac{1}{2})}
       \end{equation*}
       is a constant.
       \item Suppose $2|n$ and the local Hermitian space $\BV_v$ admits a self-dual lattice; let $\Phi'_v$ be its characteristic function. In this case, for $a\in\mathcal{O}_{F_v}$, we have
        \begin{equation*}
        W_{a,v}(s,1,\Phi'_v)=\left\{
        \begin{aligned}
           \nonumber
           &N_v^{-\frac{n}{2}}-N_v^{-(r+1)(s+n)-\frac{1}{2}} \ \ \ a\notin \Nm(E_v^\times)\cdot(-1)^\frac{n}{2}d_\BV;\\
           &N_v^{-\frac{n}{2}}+N_v^{-(r+1)(s+n)-\frac{1}{2}} \ \ \ a\in \Nm(E_v^\times)\cdot(-1)^\frac{n}{2}d_\BV,\\
        \end{aligned}
          \right.
       \end{equation*}
       and the difference
       \begin{equation*}
           W_{a,v}(s,1,\Phi'_v)- W_{a,v}(s,1,\Phi_v)=N_v^{-\frac{n}{2}}-N_v^{-(n+\frac{1}{2})}
       \end{equation*}
       is a constant.
    \end{enumerate}
\end{corollary}
\begin{proof}
    In fact, the expressions for all the cases in this corollary can be derived using the computational methods for the corresponding situations outlined in the previous Lemma \ref{Explicit nonarchimedean}. In other words, when $2\nmid n$, the method is the same as the one in Step 1; when $2|n$, the induction method in Step 4 remains valid, and the coefficients $C_m^r$ there need no change. Thus, we omit these repetitive computations.

    Meanwhile, it is worth noting that the phenomenon that the difference in Whittaker functions provided by different Schwartz functions is a constant (independent of $a$) is not a coincidence. The reason behind it can be found in \cite[Lemma 3.1.2]{Qiu}.

    When $n$ is even, the stabilization statement needed later is valid, but it requires more than local constancy alone. Fix
    \begin{equation*}
        l\in\mathcal O_{F_v}^\times,\qquad l\notin\Nm(E_v^\times).
    \end{equation*}
    Then, for every $\Phi_v\in\mathcal S(\BV_v)$, there is an integer $r_0=r_0(\Phi_v)$ such that, whenever $v(a)\ge r_0$,
    \begin{equation}\label{stabilization of local Whittaker functions}
        W_{a,v}(s,1,\Phi_v)+W_{al,v}(s,1,\Phi_v)
        =2W_{0,v}(s,1,\Phi_v)
    \end{equation}
    as meromorphic functions of $s$.

    To justify this, first take $\Phi_v=1_{\Lambda_v}$. Multiplication by $l$ exchanges the two norm classes occurring in the two cases of the explicit formula in Lemma \ref{Explicit nonarchimedean}; the terms depending on $a$ therefore cancel, and the sum is $2W_{0,v}(s,1,1_{\Lambda_v})$. The same holds for the characteristic function of any sufficiently small homothetic lattice by the covariance under $m(c)$. For general $\Phi_v$, choose such a lattice $L$ on which $\Phi_v$ is constant and write
    \begin{equation*}
        \Phi_v=\Phi_v(0)1_L+\Phi_v^0,\qquad \Phi_v^0(0)=0.
    \end{equation*}
    It remains to treat $\Phi_v^0$. Since $\Phi_v^0(0)=0$, the function
    \begin{equation*}
      b\longmapsto \delta_v(wn(b))^s r(wn(b))\Phi_v^0(0)
    \end{equation*}
    is compactly supported in $b\in F_v$: for $b$ in a fixed compact set this is clear, while for $v(b)$ sufficiently negative the matrix identity \eqref{matrix equation} reduces the value at the origin to the action of a compact element arbitrarily close to the central element $-1$, and local constancy together with $\Phi_v^0(0)=0$ makes it vanish. Its Fourier transform is therefore locally constant near $t=0$. Hence
    \begin{equation*}
      W_{t,v}(s,1,\Phi_v^0)=W_{0,v}(s,1,\Phi_v^0)
    \end{equation*}
    for all sufficiently small $t$, first in the half-plane of convergence and then as a meromorphic identity in $s$. This is the same argument as \cite[Lemma 3.1.2]{Qiu}. Applying it to $t=a$ and $t=al$, and combining it with the homothetic-lattice case, proves \eqref{stabilization of local Whittaker functions}.
\end{proof}

In fact, the stabilization identity above is the only local input from this paragraph that will be used for the exceptional $B$-families in the third paper of the series.  Its proof has a slightly stronger consequence than the formulation for Schwartz functions: once a smooth local section has been supplied on the closed Bruhat cell and has the same homogeneous Borel covariance, the identical tail-cancellation argument applies to that section.  Indeed, on the open-cell tail the proof uses only smoothness and the sign of the ramified quadratic character under the non-norm unit $l$, while on the remaining compact set it uses only local constancy of the additive characters.  Thus no Schwartz preimage for the exceptional radial coefficient function is required.  In the global comparison of the third paper, the exceptional terms are first kept inside the holomorphic automorphic residual, which supplies precisely such a smooth closed-cell completion; the stabilization identity therefore applies to its completed local two-norm symmetrization.  This does not assert that an isolated local summand is automorphic.

\subsection{Explicit computations of intertwining part}
Recall that we defined a new function \ref{c sequence} and the associated Eisenstein series \ref{C sequence}, i.e., if $\Phi_v=\phi_{1,v}\otimes\phi_{2,v}$,
\begin{equation*}
    c_{\Phi_v}(g,y)=r(g)\phi_{1,v}(y){W^\circ_{0,v}}'(0,g,\phi_{2,v})+\log\delta(g_v)r(g)\Phi_v(y),
\end{equation*}
where the normalization ${W^\circ_{0,v}}'(0,g,\phi_{2,v})$ is introduced in \ref{Normalized Whittaker function 0}. Note that this definition extends by linearity to general $\Phi$. In the final computation, we also need to use the explicit formula of $c_{\Phi_v}(g,y)$ for some fixed $g\in\UU(F_v)$. Thus, we have the following lemma.

\begin{lemma}\label{Explicit good intertwining}
    Recall our choice of the Schwartz function in Section \ref{Choice of Schwartz function}.
    \begin{enumerate}
        \item For any archimedean place $v$,
        \begin{equation*}
            c_{\Phi_v}(g,y)=0,\quad g\in\UU(\RR),\ y\in W^\perp_v.
        \end{equation*}
        \item For any non-archimedean place $v$ and $y\in W^\perp_v$, we have
        \begin{equation*}
            c_{\Phi_v}(1,y)=\Phi_v(y)\cdot\log|d_v|+\left\{
    \begin{aligned}
        \nonumber
        &\Phi_v(y)\cdot\frac{2(|d_v|-1)}{(1+N_v^{-1})(1-N_v)}\log N_v \ \ \ \mathrm{if}\ E_v/F_v\ \mathrm{inert};\\
        &0\ \ \ \mathrm{otherwise}.\\
    \end{aligned}
    \right.
        \end{equation*}
    \end{enumerate}
\end{lemma}
\begin{proof}
    Under the compatible orthogonal factorization of the local lattice and Schwartz function fixed above, the nontrivial factor is rank two, so this calculation reduces to the $n=1$ computation of \cite[Lemma 7.6]{YZ1}. Without that factorization, the phrase ``without loss of generality'' would be too strong.

    If $v$ is archimedean, it is sufficient to check that
    \begin{equation*}
        W^\circ_{0,v}(s,g,\phi_{2,v})=\delta(g)^{-s}r(g)\phi_{2,v}(0),\ g\in\UU(F_v).
    \end{equation*}
    Note that the behaviors of both sides under the left action of $P(\RR)$ and the right action of the maximal compact subgroup $K(F_v)\subset\UU(F_v)$ are the same, which implies two sides are equal up to a constant possibly depending on $s$. Then it is sufficient to check that
    $W^\circ_{0,v}(s,1,\phi_{2,v})=1$, which can be checked from the computation in Lemma \ref{Explicit archimedean}.

    Assume that $v$ is non-archimedean in the following. If $v$ is unramified in $E/F$, the computation is entirely the same as \cite[Lemma 7.6(2)]{YZ1}, as the self-dual lattice in our case corresponds to the maximal order in quaternion algebra. See Remark \ref{quaternion and Hermitian remark}.

It remains to check the case when $v$ is ramified in $E/F$. On the chosen
rank-two modular factor, apply the general ramified formula in
\cite[Lemma 7.6(2)]{YZ1} to the actual lattice-standard Schwartz function.
Its first scalar correction is
\begin{equation*}
 \Phi_v(y)\frac{|d_vq(\mathfrak j_v)|-1}{1-N_v}\log N_v,
\end{equation*}
and the remaining term $\alpha_v(y)$ is a sum indexed by
$0\leq m\leq v(d_v)-1$. Under the standing hypotheses
\begin{equation*}
 |d_v|=1,\qquad |q(\mathfrak j_v)|=1,\qquad v(D_v)=1,
\end{equation*}
the scalar correction is zero and the indexing set defining $\alpha_v$ is
empty. The complementary orthogonal factor contributes only its
characteristic function. Hence
\begin{equation*}
 c_{\Phi_v}(1,y)=0,
 \qquad y\in W_v^\perp,
\end{equation*}
which is the claimed ramified formula. No passage from an auxiliary
Schwartz function by linearity is required.

\end{proof}

Finally, Lemma \ref{Explicit archimedean k series}, \ref{Explicit nonarchimedean k series}, \ref{Explicit archimedean}, \ref{Explicit nonarchimedean} and \ref{Explicit good intertwining} in this chapter provide the explicit expressions for all pseudo-theta series and pseudo-Eisenstein series appearing in Theorems \ref{Projection of I'} and \ref{Projection of J'}, together with their corresponding coefficients including $c_0$ and $c_3$. This completes the proof of the main Theorem \ref{main theorem of derivative series} of the paper.

\

\noindent \small{School of mathematical sciences, Peking University, Beijing 100871, China}

\noindent \small{\it Email: ziqiguo0603@pku.edu.cn}

\end{document}